\documentclass[seceqn,number,citesort,dvips]{arxbj}
\usepackage{cursive}  

\aid{0}
\volume{16}
\issue{4}
\pubyear{2010}
\firstpage{1385}
\lastpage{1414}
\doi{10.3150/10-BEJ254}

\makeatletter
\newtheorem{theorem}{Theorem}[section]
\newtheorem{proposition}{Proposition}[section]
\newtheorem{lemma}[proposition]{Lemma}
\newremark{example}{Example}
\newremark{remark}{Remark}
\newremark{remarks}{Remarks}
\newtheorem{corollary}[proposition]{Corollary}
\makeatother

\begin{document}
\begin{frontmatter}

\title{Concentration of empirical distribution functions with
applications to~non-i.i.d.~models}
\runtitle{Concentration of empirical distribution functions}

\begin{aug}
\author[a]{\fnms{S.G.} \snm{Bobkov}\corref{}\thanksref{a}\ead[label=e1]{bobkov@math.umn.edu}} \and
\author[b]{\fnms{F.} \snm{G\"otze}\thanksref{b}\ead[label=e2]{goetze@mathematik.uni-bielefeld.de}}
\runauthor{S.G. Bobkov and F. G\"otze}
\address[a]{School of Mathematics, University of Minnesota, 127 Vincent
Hall, 206 Church St. S.E., Minneapolis, MN~55455, USA. \printead{e1}}
\address[b]{Fakult\"at f\"ur Mathematik, Universit\"at Bielefeld,
Postfach 100131, 33501 Bielefeld, Germany.\\ \printead{e2}}
\pdfauthor{S.G. Bobkov, F. Gotze}
\end{aug}

\received{\smonth{7} \syear{2008}}
\revised{\smonth{9} \syear{2009}}

%
\begin{abstract}
The concentration of empirical measures is studied for dependent data,
whose joint distribution satisfies Poincar\'e-type or logarithmic
Sobolev inequalities. The general concentration results are then
applied to
spectral empirical distribution functions associated with
high-dimensional random matrices.
\end{abstract}
%
%
\pdfkeywords{empirical measures, logarithmic Sobolev inequalities,
Poincare-type inequalities, random matrices, spectral distributions}

\end{frontmatter}
%
\section{Introduction}\label{intr}

Let $(X_1,\ldots,X_n)$ be a random vector in $\mathbf{R}^n$ with distribution
$\mu$. We study rates of approximation of the average marginal
distribution function
\[
F(x) = \mathbf{E}F_n(x) = \frac{1}{n} \sum_{i=1}^n \mathbf{P}\{X_i
\leq x\}
\]
by the empirical distribution function
\[
F_n(x) = \frac{1}{n} \operatorname{card}\{i \leq n\dvtx X_i \leq x\},\qquad
x \in\mathbf{R}.
\]
We shall measure the distance between $F$ and $F_n$ by means of the
(uniform) Kolmogorov metric $\|F_n - F\| = \sup_x |F_n(x) - F(x)|$,
as well as by means of the $L^1$-metric
\[
W_1(F_n,F) = \int_{-\infty}^{+\infty} |F_n(x) - F(x)| \,\mathrm{d}x.
\]
The latter, also called the Kantorovich--Rubinstein distance,
may be interpreted as the minimal cost needed to transport
the empirical measure $F_n$ to $F$ with cost function $d(x,y) = |x-y|$
(the price paid to transport the point $x$ to the point $y$).

The classical example is the case where all $X_i$'s are independent
and identically distributed (i.i.d.), that is, when $\mu$ represents
a product measure on $\mathbf{R}^n$ with equal marginals, say, $F$.
If it has no atoms, the distributions of the random variables
$T_n = \sqrt{n} \|F_n - F\|$ are weakly convergent to the Kolmogorov
law. Moreover, by the Dvoretzky--Kiefer--Wolfowitz theorem, the r.v.'s
$T_n$ are uniformly sub-Gaussian and, in particular,
$\mathbf{E}\|F_n - F\| \leq\frac{C}{\sqrt{n}}$ up to a universal
factor $C$ (\cite{D-K-W}; cf. \cite{Mas} for history and sharp bounds).
This result, together with the
related invariance principle, has a number of extensions to the
case of dependent observations, mainly in terms of mixing conditions
imposed on a stationary process; see, for example, \cite{Se,K,Y}.

On the other hand, the observations $X_1,\ldots,X_n$ may also be
generated by non-tri\-vial functions of independent random variables.
Of particular importance are random symmetric matrices
$(\frac{1}{\sqrt{n}} \xi_{jk})$, $1 \leq j,k \leq n$,
with i.i.d. entries above and on the diagonal. Arranging their
eigenvalues $X_1 \leq\cdots\leq X_n$ in increasing order,
we arrive at the spectral empirical measures $F_n$. In this case,
the mean $F = \mathbf{E}F_n$ also depends on $n$ and converges to the
semicircle law under appropriate moment assumptions on $\xi_{jk}$ (cf.,
for example, \cite{P,G-T1}).

The example of matrices strongly motivates the study of deviations of
$F_n$ from the mean $F$ under general analytical hypotheses on the
joint distribution of the observations, such as Poincar\'e or logarithmic
Sobolev inequalities. A probability measure $\mu$ on $\mathbf{R}^n$
is said
to satisfy a \textit{Poincar\'e-type} or \textit{spectral gap
inequality with constant}
$\sigma^2$ ($\sigma\geq0$) if, for any bounded smooth function $g$
on $\mathbf{R}^n$ with gradient $\nabla g$,
%
\begin{equation}\label{for1.1}
\operatorname{Var}_\mu(g) \leq\sigma^2 \int|\nabla g|^2 \,\mathrm{d}\mu.
\end{equation}
In this case, we write $\operatorname{PI}(\sigma^2)$ for short. Similarly, $\mu$
satisfies a logarithmic Sobolev inequality with constant $\sigma^2$
and we write $\operatorname{LSI}(\sigma^2)$ if, for all bounded smooth $g$,
%
\begin{equation}\label{for1.2}
\operatorname{Ent}_\mu(g^2) \leq2\sigma^2 \int|\nabla g|^2 \,\mathrm{d}\mu.
\end{equation}
Here, as usual, $\operatorname{Var}_\mu(g) = \int g^2 \,\mathrm{d}\mu- (\int g
\,\mathrm{d}\mu)^2$
stands for the variance of $g$ and
$\operatorname{Ent}_\mu(g) = \int g \log g \,\mathrm{d}\mu- \int g \,\mathrm{d}\mu \log
\int g \,\mathrm{d}\mu$
denotes the entropy of $g \geq0$ under the measure $\mu$.
It is well known that LSI($\sigma^2$) implies PI($\sigma^2$).

These hypotheses are crucial in the study of concentration of the
spectral empirical distributions, especially of the linear functionals
$\int f\, \mathrm{d}F_n$ with individual smooth $f$ on the line; see, for
example, the results by Guionnet and Zeitouni \cite{G-Z},
Chatterjee and Bose \cite{C-B}, Davidson and Szarek \cite{D-S} and
Ledoux \cite{L2}. A remarkable feature of this approach to spectral
analysis is that no specific knowledge about the non-explicit
mapping from a random matrix to its spectral empirical measure
is required. Instead, one may use general Lipschitz properties only,
which are satisfied by this mapping.
As for the general (not necessarily matrix) scheme, we shall only
require the hypotheses (\ref{for1.1}) or (\ref{for1.2}). In particular, we derive the
following from (\ref{for1.1}).

\begin{theorem}\label{teo1.1} Under $\operatorname{PI}(\sigma^2)$ on $\mathbf{R}^n$ $(n
\geq2)$,
%
\begin{equation}\label{for1.3}
\mathbf{E}\int_{-\infty}^{+\infty} |F_n(x) - F(x)| \,\mathrm{d}x \leq
C\sigma\biggl(\frac{A + \log n}{n}\biggr)^{ 1/3},
\end{equation}
where $A = \frac{1}{\sigma} \max_{i,j} | \mathbf{E}X_i - \mathbf
{E}X_j|$ and
$C$ is an absolute constant.
\end{theorem}

Note that the Poincar\'e-type inequality (\ref{for1.1}) is invariant under shifts
of the measure $\mu$, while the left-hand side of (\ref{for1.3}) is not. This is
why the bound on the right-hand side of (\ref{for1.3}) should also depend on the
means of the observations.

In terms of the ordered statistics $X_1^* \leq\cdots\leq X_n^*$ of the
random vector $(X_1,\ldots,X_n)$, there is a general two-sided estimate
for the mean of the Kantorovich--Rubinstein distance:
%
\begin{equation}\label{for1.4}
\frac{1}{2n} \sum_{i=1} \mathbf{E}|X_i^* - \mathbf{E}X_i^*| \leq
\mathbf{E}W_1(F_n,F)
\leq\frac{2}{n} \sum_{i=1} \mathbf{E}|X_i^* - \mathbf{E}X_i^*|
\end{equation}
(see remarks at the end of Section \ref{sec4}).
Hence, under the conditions of Theorem~\ref{teo1.1}, one may control the local
fluctuations of $X_i^*$ (on average), which typically deviate from their
mean by not more than $C\sigma(\frac{A + \log n}{n})^{1/3}$.

Under a stronger hypothesis, such as (\ref{for1.2}), one can obtain more information
about the fluctuations of $F_n(x) - F(x)$ for individual points $x$ and
thus get some control of the Kolmogorov distance. Similarly to the
bound (\ref{for1.3}), such fluctuations will, on average, be shown to be at most
\[
\beta= \frac{(M\sigma)^{2/3}}{n^{1/3}},
\]
in the sense that $\mathbf{E}|F_n(x) - F(x)| \leq C \beta$, where $M$
is the
Lipschitz seminorm of $F$ (see Proposition \ref{pro6.3}). As for the Kolmogorov
distance, we prove the following theorem.

\begin{theorem}\label{teo1.2} Assume that $F$ has a density, bounded by
a number $M$.
Under $\operatorname{LSI}(\sigma^2)$, for any $r>0$,
%
\begin{equation}\label{for1.5}
\mathbf{P}\{\|F_n - F\| \geq r \} \leq
\frac{4}{r} \mathrm{e}^{-c (r/\beta)^3}.
\end{equation}
In particular,
%
\begin{equation}\label{for1.6}
\mathbf{E}\|F_n - F\| \leq C \beta
\log^{1/3} \biggl(1 + \frac{1}{\beta}\biggr),
\end{equation}
where $c$ and $C$ are positive absolute constants.
\end{theorem}

In both cases, the stated bounds are of order $n^{-1/3}$ up to a $\log n$
term with respect to the dimension $n$. Thus, they are not as sharp as
in the classical i.i.d. case. Indeed, our assumptions are much weaker
and may naturally lead to weaker conclusions.
Let us look at two examples illustrating the
bounds obtained in the cases that essentially differ from the i.i.d. case.

\begin{example}\label{ex1} Let $X_i$ be independent and uniformly distributed
in the intervals $(i-1,i)$, $i = 1,\ldots,n$. Their joint distribution
is a product measure, satisfying (\ref{for1.1}) and (\ref{for1.2}) with some absolute
constant $\sigma$. Clearly, $F$ is the uniform distribution in $(0,n)$
so $M = \frac{1}{n}$ and $\beta$ is of order $\frac{1}{n}$.
As is easy to see, $\mathbf{E}\|F_n - F\|$ is also of order $\frac{1}{n}$,
that is, the bound (\ref{for1.6}) is sharp up to a $\log^{1/3} n$ term.
Also, since $A$ is of order $n$, both sides of (\ref{for1.3}) are of order 1.
In particular, this shows that the quantity $A$ cannot be removed from (\ref{for1.3}).
\end{example}

\begin{example}\label{ex2} Let all $X_i = \xi$, where $\xi$ is uniformly
distributed in $[-1,1]$. Note that all random variables are identically
distributed with $\mathbf{E}X_i = 0$. The joint distribution $\mu$ represents
a uniform distribution on the main diagonal of the cube $[-1,1]^n$,
so it satisfies (\ref{for1.1}) and (\ref{for1.2}) with $\sigma= c\sqrt{n}$, where
$c$ is absolute. In this case, $F$ is a uniform distribution on $[-1,1]$,
so $M = 1/2$ and $\beta$ is of order 1. Hence, both sides of (\ref{for1.6})
are of order 1.
\end{example}

Next, we restrict the above statements to the empirical spectral
measures $F_n$ of the $n$ eigenvalues $X_1 \leq\cdots\leq X_n$
of a random symmetric matrix
$(\frac{1}{\sqrt{n}} \xi_{jk})$, ${1 \leq j,k \leq n}$, with
independent entries above and on the diagonal ($n \geq2$). Assume that
$\mathbf{E}\xi_{jk} = 0$ and ${\operatorname{Var}(\xi_{jk}) = 1}$
so that the means
$F = \mathbf{E}F_n$ converge to the semicircle law $G$ with mean zero and
variance one. The boundedness of moments of $\xi_{jk}$ of any order
will be guaranteed by (\ref{for1.1}).

\begin{theorem}\label{teo1.3} If the distributions of the $\xi_{jk}$'s
satisfy the
Poincar\'e-type inequality $\operatorname{PI}(\sigma^2)$ on the real line, then
%
\begin{equation}\label{for1.7}
\mathbf{E}\int_{-\infty}^{+\infty} |F_n(x) - F(x)| \,\mathrm{d}x \leq
\frac{C\sigma}{n^{2/3}},
\end{equation}
where $C$ is an absolute constant. Moreover, under $\operatorname{LSI}(\sigma^2)$,
%
\begin{equation}\label{for1.8}
\mathbf{E}\|F_n - G\| \leq C \biggl(\frac{\sigma}{n}\biggr)^{ 2/3}
\log^{1/3} n + \|F-G\|.
\end{equation}
\end{theorem}

By the convexity of the distance, we always have
$\mathbf{E}\|F_n - G\| \geq\|F-G\|$. In some random matrix models, the
Kolmogorov distance $\|F-G\|$ is known to tend to zero at rate at most
$n^{-2/3 + \varepsilon }$. For instance, it is true when the distributions
of the $\xi_{j,k}$'s have a non-trivial Gaussian component (see \cite{G-T-T}).
Hence, if, additionally, $\operatorname{LSI}(\sigma^2)$ is satisfied, then we
get that
for any $\varepsilon > 0$,
\[
\mathbf{E}\|F_n - G\| \leq C_{\varepsilon,  \sigma} n^{-2/3 +
\varepsilon }.
\]
It is unknown whether this bound is optimal. Note, however, that in the case
of Gaussian $\xi_{j,k}$, the distance $\|F-G\|$ is known to be of order
$1/n$ \cite{G-T2}. Therefore,
%
\begin{equation}\label{for1.9}
\mathbf{E}\|F_n - G\| \leq C \frac{\log^{1/3} n}{n^{2/3}},
\end{equation}
which is a slight improvement of a bound obtained in \cite{Ti}. In fact,
as was recently shown in \cite{B-G-T}, we have $\|F-G\| \leq C_\sigma n^{-2/3}$
in the presence of the
$\operatorname{PI}(\sigma^2)$-hypothesis. Hence, the bound (\ref{for1.9}) always holds
under $\operatorname{LSI}(\sigma^2)$ with constants depending only on $\sigma$.

It seems natural to try to relax the $\operatorname{LSI}(\sigma^2)$-hypothesis
in (\ref{for1.8}) and (\ref{for1.9}) to $\operatorname{PI}(\sigma^2)$. In this context, let us mention
a result of Chatterjee and Bose \cite{C-B}, who used Fourier transforms
to derive from $\operatorname{PI}(\sigma^2)$ a similar bound,
\[
\mathbf{E}\|F_n - G\| \leq\frac{C\sigma^{1/4}}{n^{1/2}} + 2 \|F-G\|.
\]
As for (\ref{for1.7}), let us return to the two-sided bound (\ref{for1.4}) which holds
with $X_i^* = X_i$ by the convention that the eigenvalues are listed
in increasing order. The asymptotic behavior of distributions of
$X_i$ with fixed or varying indices has been studied by many authors, especially
in the standard Gaussian case. In particular, if $i$ is fixed, while
$n$ grows, $n^{2/3} (X_i -\mathbf{E}X_i)$ converges in distribution to
(a variant of) the Tracy--Widom law so the $\mathbf{E}|X_i - \mathbf
{E}X_i|$ are of
order $n^{-2/3}$. This property still holds when $\xi_{jk}$ are
symmetric and have sub-Gaussian tails; see \cite{So} and \cite{L3} for the history
and related results. Although this rate is consistent with the bound
(\ref{for1.7}), the main contribution in the normalized sum (\ref{for1.4}) is due to
the intermediate terms (in the bulk) and their rate might be different.
It was shown by Gustavsson \cite{Gu} for the GUE model that if
$\frac{i}{n} \rightarrow t \in(0,1)$,\vspace*{-2pt} then $X_i$ is asymptotically normal
with variance of order $\frac{C(t)\log n}{n^2}$. Hence, it is not surprising
that\vspace*{-1pt} $\mathbf{E}W_1(F_n,F) \leq\frac{C (\log n)^{1/2}}{n}$, see \cite{Ti}.

The paper is organized as follows. In Section \ref{sec2}, we collect a few
direct applications of the Poincar\'e-type inequality to linear
functionals of empirical measures. They are used in Section \ref{sec3} to
complete the proof of Theorem~\ref{teo1.1}. In the next section, we discuss
deviations of $W_1(F_n,F)$ from its mean. In Section \ref{sec5}, we turn to
logarithmic Sobolev inequalities. Here, we shall adapt
infimum-convolution operators to empirical measures and apply a result
of \cite{B-G-L} on the relationship between infimum-convolution and
log-Sobolev inequalities. In Section \ref{sec6}, we illustrate this approach
in the problem of dispersion of the values of the empirical distribution
functions at a fixed point.
In Section \ref{sec7}, we derive bounds on the uniform distance similar
to (\ref{for1.5}) and (\ref{for1.6}) and give a somewhat more general form of Theorem \ref{teo1.2}.
In Section \ref{sec8}, we apply the previous results to high-dimensional
random matrices to prove Theorem \ref{teo1.3} and obtain some refinements.
Finally, since the hypotheses (\ref{for1.1}) and (\ref{for1.2}) play a crucial role
in this investigation, we collect in the Appendix a few results on
sufficient conditions for a measure to satisfy PI and LSI.


\section{Empirical Poincar\'e inequalities}\label{sec2}

We assume that the random variables $X_1,\ldots,X_n$ have a joint
distribution $\mu$ on $\mathbf{R}^n$, satisfying the Poincar\'e-type
inequality (\ref{for1.1}). For a bounded smooth function $f$ on the real line,
we apply it to
%
\begin{equation}\label{for2.1}
g(x_1,\ldots,x_n) = \frac{f(x_1) + \cdots+ f(x_n)}{n} =
\int f \,\mathrm{d}F_n,
\end{equation}
where $F_n$ is the empirical measure, defined for the `observations'
$X_1 = x_1,\ldots,X_n = x_n$. Since
%
\begin{equation}\label{for2.2}
|\nabla g(x_1,\ldots,x_n)|^2 =
\frac{f'(x_1)^2 + \cdots+ f'(x_n)^2}{n^2} = \frac{1}{n} \int f'^2 \,\mathrm{d}F_n,
\end{equation}
we obtain an integro-differential inequality, which may viewed as
an empirical Poincar\'e-type inequality for the measure $\mu$.

\begin{proposition}\label{pro2.1} Under $\operatorname{PI}(\sigma^2)$, for
any smooth
$F$-integrable function $f$ on $\mathbf{R}$ such that $f'$ belongs to
$L^2(\mathbf{R},\mathrm{d}F)$, we have
%
\begin{equation}\label{for2.3}
\mathbf{E}\bigg|\int f\,\mathrm{d}F_n - \int f\,\mathrm{d}F \bigg|^2 \leq
\frac{\sigma^2}{n} \int f'^2\,\mathrm{d}F.
\end{equation}
\end{proposition}

Recall that $F = \mathbf{E}F_n$ denotes the mean of the empirical measures.
The inequality continues to hold for all locally Lipschitz functions with
the modulus of the derivative, understood in the generalized sense,
that is,
$|f'(x)| = \limsup_{y \rightarrow x} \frac{|f(x)-f(y)|}{|x-y|}$.
As long as $\int f'^2\,\mathrm{d}F$ is finite, $\int f^2\,\mathrm{d}F$ is
also finite and (\ref{for2.3}) holds.

The latter may be extended to all $L^p$-spaces by
applying the following general lemma.

\begin{lemma}\label{lem2.2} Under $\operatorname{PI}(\sigma^2)$, any Lipschitz
function $g$ on $\mathbf{R}^n$ has a finite exponential moment: if
$\int g \,\mathrm{d}\mu= 0$ and $\|g\|_{\mathrm{Lip}} \leq1$, then
%
\begin{equation}\label{for2.4}
\int \mathrm{e}^{tg/\sigma} \,\mathrm{d}\mu\leq\frac{2 + t}{2 - t},\qquad
0 < t < 2.
\end{equation}
Moreover, for any locally Lipschitz $g$ on $\mathbf{R}^n$ with $\mu
$-mean zero,
%
\begin{equation}\label{for2.5}
\|g\|_p \leq\sigma p \|\nabla g\|_p,\qquad  p \geq2.
\end{equation}
\end{lemma}

More precisely, if $|\nabla g|$ is in $L^p(\mu)$, then so is $g$ and
(\ref{for2.5}) holds true with the standard notation
$\|g\|_p = (\int|g|^p \,\mathrm{d}\mu)^{1/p}$ and
$\|\nabla g\|_p = (\int|\nabla g|^p \,\mathrm{d}\mu)^{1/p}$
for $L^p(\mu)$-norms. The property of being locally Lipschitz
means that the function $g$ has a finite Lipschitz seminorm
on every compact subset of $\mathbf{R}^n$.

In the concentration context, a variant of the first part of the lemma
was first established by Gromov and Milman in \cite{G-M} and
independently in dimension $1$ by Borovkov and Utev \cite{B-U}.
Here, we follow \cite{B-L1}, Proposition 4.1,
to state (\ref{for2.4}). The second inequality, (\ref{for2.5}), may be derived by similar
arguments; see also \cite{B-L2}, Theorem 4.1, for an extension to the case
of Poincar\'e-type inequalities with weight.

Now, for functions $g = \int f\,\mathrm{d}F_n$ as in (\ref{for2.1}), in view of (\ref{for2.2}),
we may write
\[
|\nabla g|^p = \frac{1}{n^{p/2}} \biggl(\int f'^2\,\mathrm{d}F_n\biggr)^{ p/2}
\leq\frac{1}{n^{p/2}} \int|f'|^p\,\mathrm{d}F_n
\]
so that $\mathbf{E}_\mu|\nabla g|^p \leq\frac{1}{n^{p/2}} \int
|f'|^p\,\mathrm{d}F$.
Applying (\ref{for2.5}) and (\ref{for2.4}) with $t=1$, we obtain the following
proposition.

\begin{proposition}\label{pro2.3} Under $\operatorname{PI}(\sigma^2)$, for
any smooth function $f$ on $\mathbf{R}$ such that $f'$ belongs to
$L^p(\mathbf{R},\mathrm{d}F)$, $p \geq2$,
\[
\mathbf{E}\biggl|\int f\,\mathrm{d}F_n - \int f\,\mathrm{d}F \biggr|^p \leq
\frac{(\sigma p)^p}{n^{p/2}} \int|f'|^p\,\mathrm{d}F.
\]
In addition, if $|f'| \leq1$, for all $h>0$,
\[
\mu\biggl\{\biggl|\int f\,\mathrm{d}F_n - \int f\,\mathrm{d}F \biggr| \geq h\biggr\}
\leq6 \mathrm{e}^{-nh/\sigma}.
\]
\end{proposition}

The empirical Poincar\'e-type inequality (\ref{for2.3}) can be rewritten
equivalently if we integrate by parts the first integral as
$
\int f\,\mathrm{d}F_n - \int f\,\mathrm{d}F = -\int f'(x) (F_n(x) - F(x)) \,\mathrm{d}x.
$
At this step, it is safe to assume that $f$ is continuously
differentiable and is constant near $-\infty$ and $+\infty$.
Replacing $f'$ with $f$, we arrive at
%
\begin{equation}\label{for2.6}
\mathbf{E}\biggl| \int f(x) \bigl(F_n(x) - F(x)\bigr) \,\mathrm{d}x \biggr|^{ 2} \leq
\frac{\sigma^2}{n} \int f^2\,\mathrm{d}F
\end{equation}
for any continuous, compactly supported function $f$ on the line.
In other words, the integral operator
$Kf(x) = \int_{-\infty}^{+\infty} K(x,y) f(y)\, \mathrm{d}y$
with a (positive definite) kernel
\[
K(x,y) = \mathbf{E}\bigl(F_n(x) - F(x)\bigr)\bigl(F_n(y) - F(y)\bigr) = \operatorname{cov}(F_n(x),F_n(y))
\]
is continuous and defined on a dense subset of
$L^2(\mathbf{R},\mathrm{d}F(x))$, taking values in $L^2(\mathbf{R},\mathrm{d}x)$. It
has the operator
norm $\|K\| \leq\frac{\sigma}{\sqrt{n}}$, so it may be continuously
extended to the space $L^2(\mathbf{R},\mathrm{d}F)$ without a change of the norm.
In the following, we will use a particular case of (\ref{for2.6}).

\begin{corollary}\label{cor2.4} Under $\operatorname{PI}(\sigma^2)$, whenever $a<b$,
we have
\[
\mathbf{E}\biggl|\int_a^b \bigl(F_n(x) - F(x)\bigr) \,\mathrm{d}x \biggr| \leq
\frac{\sigma}{\sqrt{n}} \sqrt{F(b)-F(a)}.
\]
\end{corollary}


\section[Proof of Theorem 1.1]{Proof of Theorem \protect\ref{teo1.1}}\label{sec3}

We shall now study the concentration properties of empirical measures
$F_n$ around their mean $F$ based on Poincar\'e-type inequalities.
In particular, we shall prove Theorem~\ref{teo1.1}, which provides a bound
on the mean of the Kantorovich--Rubinstein distance
\[
W_1(F_n,F) = \int_{-\infty}^{+\infty} |F_n(x) - F(x)| \,\mathrm{d}x.
\]
Note that it is homogeneous of order 1 with respect to the random vector
$(X_1,\ldots,X_n)$.

We first need a general observation.

\begin{lemma}\label{lem3.1} Given distribution functions $F$ and $G$,
for all real $a<b$ and a natural number $N$,
\[
\int_a^b |F(x) - G(x)| \,\mathrm{d}x \leq\sum_{k=1}^N
\biggl|\int_{a_{k-1}}^{a_k} \bigl(F(x) - G(x)\bigr) \,\mathrm{d}x\biggr| +
\frac{2 (b-a)}{N},
\]
where $a_k = a + (b-a )\frac{k}{N}$.
\end{lemma}

\begin{pf} Let $I$ denote the collection of those indices $k$ such
that in the $k$th subinterval $\Delta_k = (a_{k-1},a_k)$, the function
$\varphi(x) = F(x) - G(x)$ does not change sign. Let $J$ denote
the collection of the remaining indices. Then, for $k \in I$,
\[
\int_{\Delta_k} |F(x) - G(x)| \,\mathrm{d}x =
\biggl|\int_{\Delta_k} \bigl(F(x) - G(x)\bigr) \,\mathrm{d}x\biggr|.
\]
In the other case $k \in J$, since $\varphi$ changes sign on $\Delta_k$,
we may write
\begin{eqnarray*}
\sup_{x \in\Delta_k} |\varphi(x)| & \leq&
\operatorname{Osc}_{\Delta_k}(\varphi) \equiv
\sup_{x,y \in\Delta_k} \bigl(\varphi(x) - \varphi(y)\bigr) \\
& \leq& \operatorname{Osc}_{\Delta_k}(F) + \operatorname{Osc}_{\Delta_k}(G) =
F(\Delta_k) + G(\Delta_k),
\end{eqnarray*}
where, in the last step, $F$ and $G$ are treated as probability measures.
Hence, in this case,
$\int_{\Delta_k} |F(x) - G(x)| \,\mathrm{d}x \leq
(F(\Delta_k) + G(\Delta_k)) |\Delta_k|$.
Combining the two bounds and using $|\Delta_k| = \frac{b-a}{N}$,
we get that
\begin{eqnarray*}
&&\int_a^b |F(x) - G(x)| \,\mathrm{d}x \\
&&\quad \leq\sum_{k \in I} \biggl|\int_{\Delta_k} \bigl(F(x) - G(x)\bigr) \,\mathrm{d}x\biggr| +
\sum_{k \in J} \bigl(F(\Delta_k) + G(\Delta_k)\bigr) |\Delta_k| \\
&&\quad \leq
\sum_{k=1}^N \biggl|\int_{\Delta_k} \bigl(F(x) - G(x)\bigr) \,\mathrm{d}x\biggr| +
\frac{b-a}{N} \sum_{k=1}^N \bigl(F(\Delta_k) + G(\Delta_k)\bigr).
\end{eqnarray*}
\upqed\end{pf}

\begin{remark*} As the proof shows, the lemma may be
extended to an arbitrary partition
$a = a_0 < a_1 < \cdots< a_N = b,$ as follows:
\[
\int_a^b |F(x) - G(x)| \,\mathrm{d}x \leq\sum_{k=1}^N
\biggl|\int_{a_{k-1}}^{a_k} \bigl(F(x) - G(x)\bigr) \,\mathrm{d}x\biggr| +
2 \max_{1 \leq k \leq N} (a_k - a_{k-1}).
\]

Let us now apply the lemma to the space $(\mathbf{R}^n,\mu)$ satisfying
a Poincar\'e-type inequality. Consider the partition of the interval
$[a,b]$ with $\Delta_k = (a_{k-1},a_k)$, as in Lemma \ref{lem3.1}.
By Corollary~\ref{cor2.4},
\begin{eqnarray*}
\mathbf{E}\int_a^b |F_n(x) - F(x)| \,\mathrm{d}x & \leq& \sum_{k=1}^N
\mathbf{E}\biggl|\int_{\Delta_k} \bigl(F_n(x) - F(x)\bigr) \,\mathrm{d}x\biggr| +
\frac{2(b-a)}{N} \\
& \leq&
\frac{\sigma}{\sqrt{n}}\sum_{k=1}^N \sqrt{F(\Delta_k)} +
\frac{2(b-a)}{N}.
\end{eqnarray*}
By Cauchy's inequality,
$\sum_{k=1}^N \sqrt{F(\Delta_k)} \leq\sqrt{N}
(\sum_{k=1}^N F(\Delta_k))^{1/2} \leq\sqrt{N}$, hence,
\[
\mathbf{E}\int_a^b |F_n(x) - F(x)| \,\mathrm{d}x \leq
\frac{\sigma\sqrt{N}}{\sqrt{n}} + \frac{2(b-a)}{N}.
\]

Now, let us rewrite the right-hand side as
$\frac{\sigma}{\sqrt{n}} (\sqrt{N} + \frac{c}{N})$
with parameter $c = \frac{2(b-a)}{\sigma/\sqrt{n}}$ and
optimize it over $N$. On the half-axis $x > 0$, introduce the function
$\psi(x) = \sqrt{x} + \frac{c}{x}$ ($c>0$). It has derivative
$\psi'(x) = \frac{1}{2\sqrt{x}} - \frac{c}{x^2}$, therefore
$\psi$ is decreasing on $(0,x_0]$ and is increasing on $[x_0,+\infty)$,
where $x_0 = (2c)^{2/3}$. Hence, if $c \leq\frac{1}{2}$,
we have
\[
\inf_{N} \psi(N) = \psi(1) = 1 + c \leq1 + c^{1/3}.
\]
If $c \geq\frac{1}{2}$, then the argmin lies in $[1,+\infty)$.
Choose $N = [x_0] + 1 = [(2c)^{2/3}] + 1$ so that $N \geq2$
and $N-1 \leq x_0 < N \leq x_0 + 1$. Hence, we get
\[
\psi(N) \leq\bigl(\sqrt{x_0} + 1\bigr) + \frac{c}{x_0} =
1 + \psi(x_0) = 1 + \frac{3}{2^{2/3}} c^{1/3}.
\]
Thus, in both cases,
$\inf_{N} \psi(N) \leq1 + \frac{3}{2^{2/3}} c^{1/3}
\leq1 + 3 (\frac{b-a}{\sigma/\sqrt{n}})^{ 1/3}$
and we arrive at the following corollary.
\end{remark*}

\begin{corollary}\label{cor3.2} Under $\operatorname{PI}(\sigma^2)$,
for all $a<b$,
\[
\mathbf{E}\int_a^b |F_n(x) - F(x)| \,\mathrm{d}x \leq \frac{\sigma}{\sqrt{n}} \biggl[
1 + 3 \biggl(\frac{b-a}{\sigma/\sqrt{n}}\biggr)^{ 1/3}\biggr].
\]
\end{corollary}

The next step is to extend the above inequality to the whole real line.
Here, we shall use the exponential integrability of the measure $F$.

\begin{pf*}{Proof of Theorem~\ref{teo1.1}}
Recall that the measure $\mu$ is controlled by using two independent
parameters: the constant $\sigma^2$ and $A$, defined by
\[
|\mathbf{E}X_i - \mathbf{E}X_j| \leq A\sigma,\qquad  1 \leq i,j \leq n.
\]
One may assume, without loss of generality, that
$-A\sigma\leq\mathbf{E}X_i \leq A\sigma$ for all $i \leq n$.

Lemma 2.2 with $g(x) = x_i$, $t=1$ and Chebyshev's inequality give,
for all $h>0$,
\[
\mathbf{P}\{X_i - \mathbf{E}X_i \geq h\} \leq3 \mathrm{e}^{-h/\sigma},\qquad
\mathbf{P}\{X_i - \mathbf{E}X_i \leq-h\} \leq3 \mathrm{e}^{-h/\sigma}.
\]
Therefore, whenever $h \geq A\sigma$,
\[
\mathbf{P}\{X_i \geq h\} \leq3 \mathrm{e}^{-(h - A\sigma)/\sigma},\qquad
\mathbf{P}\{X_i \leq-h\} \leq3 \mathrm{e}^{-(h- A\sigma)/\sigma}.
\]
Averaging over all $i$'s, we obtain similar bounds for the measure
$F$, that is, $1 - F(h) \leq3 \mathrm{e}^{-(h - A\sigma)/\sigma}$ and
$F(-h) \leq3 \mathrm{e}^{-(h- A\sigma)/\sigma}$.
After integration, we get
\[
\int_h^{+\infty} \bigl(1-F(x)\bigr) \,\mathrm{d}x \leq3\sigma \mathrm{e}^{-(h - A\sigma)/\sigma},
\qquad
\int_{-\infty}^{-h} F(x) \,\mathrm{d}x \leq3\sigma \mathrm{e}^{-(h - A\sigma)/\sigma}.
\]
Using $|F_n(x) - F(x)| \leq(1-F_n(x)) + (1-F(x))$ so that
$\mathbf{E}|F_n(x) - F(x)| \leq2(1-F(x))$, we get that
\[
\mathbf{E}\int_h^{+\infty} |F_n(x) - F(x)| \,\mathrm{d}x \leq
6\sigma \mathrm{e}^{-(h - A\sigma)/\sigma}
\]
and similarly for the half-axis $(-\infty,-h)$.
Combining this bound with Corollary \ref{cor3.2}, with $[a,b] = [-h,h]$,
we obtain that, for all $h \geq A\sigma$,
\[
\mathbf{E}\int_{-\infty}^{+\infty} |F_n(x) - F(x)| \,\mathrm{d}x \leq
\frac{\sigma}{\sqrt{n}} \biggl[
1 + 6 \biggl(\frac{h}{\sigma/\sqrt{n}}\biggr)^{ 1/3}\biggr] +
12 \sigma \mathrm{e}^{-(h - A\sigma)/\sigma}.
\]
Substituting $h = (A+t)\sigma$ with arbitrary $t \geq0$, we get that
\[
\mathbf{E}\int_{-\infty}^{+\infty} |F_n(x) - F(x)| \,\mathrm{d}x \leq
\frac{\sigma}{\sqrt{n}} \biggl[
1 + 6 \bigl((A+t)\sqrt{n} \bigr)^{1/3} + 12\sqrt{n} \mathrm{e}^{-t}\biggr].
\]
Finally, the choice $t = \log n$ leads to the desired estimate
\[
\mathbf{E}\int_{-\infty}^{+\infty} |F_n(x) - F(x)| \,\mathrm{d}x \leq
C\sigma\biggl(\frac{A + \log n}{n}\biggr)^{ 1/3}.
\]
\upqed\end{pf*}


\section{Large deviations above the mean}\label{sec4}
\setcounter{equation}{0}

In addition to the upper bound on the mean of the Kantorovich--Rubinstein
distance $W_1(F_n,F)$, one may wonder how to bound large
deviations of $W_1(F_n,F)$ above the mean. To this end, the
following general observation may be helpful.

\begin{lemma}\label{lem4.1} For all points $x = (x_1,\ldots,x_n)$,
$x' = (x_1',\ldots,x_n')$ in $\mathbf{R}^n$,
we have
\[
W_1(F_n,F_n') \leq\frac{1}{\sqrt{n}} \|x-x'\|,
\]
where
$F_n = \frac{\delta_{x_1} + \cdots+ \delta_{x_n}}{n}$,
$F_n' = \frac{\delta_{x_1'} + \cdots+ \delta_{x_n'}}{n}$.
\end{lemma}

In other words, the canonical map $T$ from $\mathbf{R}^n$ to the space of
all probability measures on the line, which assigns to each point
an associated empirical measure, has a Lipschitz seminorm
${\leq}\frac{1}{\sqrt{n}}$ with respect to the Kantorovich--Rubinstein
distance. As usual, the Euclidean space $\mathbf{R}^n$ is equipped
with the
Euclidean metric
\[
\|x-x'\| = \sqrt{|x_1 - x_1'|^2 + \cdots+ |x_n - x_n'|^2}.
\]

Denote by $Z_1$ the collection of all (Borel) probability measures
on the real line with finite first moment. The Kantorovich--Rubinstein
distance in $Z_1$ may equivalently be defined (cf. \cite{V,D}) by
\[
W_1(G,G') = \inf_\pi\int|u-u'|\, \mathrm{d}\pi(u,u'),
\]
where the infimum is taken over all (Borel) probability measures $\pi$
on $\mathbf{R}\times\mathbf{R}$ with marginal distributions $G$ and
$G'$. In case of
empirical measures $G = F_n$, $G' = F_n'$, associated to the points
$x,x' \in\mathbf{R}^n$, let $\pi_0$ be the discrete measure on the
pairs $(x_i,x_i')$, $1 \leq i \leq n$, with point masses $\frac{1}{n}$.
Therefore, by Cauchy's inequality,
\[
W_1(Tx,Tx') \leq\int|u-u'|\, \mathrm{d}\pi_0(u,u') =
\frac{1}{n} \sum_{i=1} |x_i - x_i'|
\leq
\frac{1}{\sqrt{n}} \biggl(\sum_{i=1} |x_i - x_i'|^2\biggr)^{ 1/2}.
\]
This proves Lemma \ref{lem4.1}.

Thus, the map $T\dvtx\mathbf{R}^n \rightarrow Z_1$ has the Lipschitz seminorm
$\|T\|_{\mathrm{Lip}} \leq\frac{1}{\sqrt{n}}$. As a consequence, given a
probability measure $\mu$ on $\mathbf{R}^n$, this map transports many
potential
properties of $\mu$, such as concentration, to the space $Z_1$, equipped
with the Borel probability measure $\Lambda= \mu T^{-1}$. Note that it is
supported on the set of all probability measures with at most $n$ atoms.
In particular, if $\mu$ satisfies a concentration inequality of the form
%
\begin{equation}\label{for4.1}
1 - \mu(A^h) \leq\alpha(h),\qquad   h > 0,
\end{equation}
in the class of all Borel sets $A$ in $\mathbf{R}^n$ with measure
$\mu(A) \geq\frac{1}{2}$ (where $A^h$ denotes an open Euclidean
$h$-neighborhood of $A$), then $\Lambda$ satisfies a similar
(and in fact stronger) property
\[
1 - \Lambda(B^h) \leq\alpha\bigl(h\sqrt{n} \bigr),\qquad   h > 0,
\]
in the class of all Borel sets $B$ in $Z_1$ with measure
$\Lambda(B) \geq\frac{1}{2}$ (with respect to the
Kantorovich--Rubinstein distance). In other words, an optimal so-called
concentration function $\alpha= \alpha_\mu$ in (\ref{for4.1}) for the measure
$\mu$ is related to the concentration function of $\Lambda$ by
%
\begin{equation}\label{for4.2}
\alpha_\Lambda(h) \leq\alpha_\mu\bigl(h\sqrt{n} \bigr),\qquad   h > 0.
\end{equation}

Now, in general, the concentration function has a simple functional
description as
\[
\alpha_\mu(h) = \sup\mu\{g - m(g) \geq h\},
\]
where the sup is taken over all Lipschitz functions $g$ on $\mathbf{R}^n$
with $\|g\|_{\mathrm{Lip}} \leq1$ and where $m(g)$ stands for a median
of $g$ under $\mu$. (Actually, this holds for abstract metric spaces.)
The concentration function may therefore be controlled
by Poincar\'e-type inequalities in terms of $\sigma^2$
(the Gromov--Milman theorem). Indeed, since the quantity $g - m(g)$
is translation invariant, one may assume that $g$ has mean zero.
By Lemma \ref{lem2.2} with $t=1$, we get
$\mu\{g \leq-\sigma h\} \leq3\mathrm{e}^{-h} < \frac{1}{2}$, provided that
$h > \log6$, which means that any median of $g$ satisfies
$m(g) \geq-\sigma\log6$. Therefore, again by Lemma \ref{lem2.2}, for any
$h > \log6$,
\[
\mu\{g - m(g) \geq\sigma h\} \leq
\mu\{g \geq\sigma(h - \log6)\} \leq3 \cdot6 \cdot \mathrm{e}^{-h}
\]
so that
%
\begin{equation}\label{for4.3}
\alpha_\mu(\sigma h) \leq18 \mathrm{e}^{-h}.
\end{equation}
The latter also automatically holds in the interval $0 \leq h \leq\log
6$. In fact, by a more careful application of the Poincar\'e-type
inequality, the concentration bound (\ref{for4.3}) may be further improved to
$\alpha_\mu(\sigma h) \leq C \mathrm{e}^{-2h}$ (see \cite{Bo2}), but this is not
crucial for our purposes.

Thus, combining (\ref{for4.2}) with (\ref{for4.3}), we may conclude that under
PI($\sigma^2$),
\[
\alpha_{\Lambda}(h) \leq18 \mathrm{e}^{-h\sqrt{n}/\sigma},\qquad  h>0.
\]

Now, in the setting of Theorem~\ref{teo1.1}, consider on $Z_1$ the distance function
$g(H) = W_1(H,F)$. It is Lipschitz (with Lipschitz seminorm $1$) and has
the mean
$\mathbf{E}_\Lambda g = \mathbf{E}_\mu W_1(F_n,F) \leq a$,
where $a = C\sigma(\frac{A + \log n}{n})^{1/3}$.
Hence, $m(g) \leq2a$ under the measure $\Lambda$ and
for any $h>0$,
\[
\Lambda\{g \geq2a + h\} \leq\Lambda\{g - m(g) \geq h\}
\leq\alpha_{\Lambda}(h) \leq18 \mathrm{e}^{-h\sqrt{n}/\sigma}.
\]
We can summarize as follows.

\begin{proposition}\label{pro4.2} If a random vector $(X_1,\ldots
,X_n)$ in
$\mathbf{R}^n$, $n \geq2$, has distribution satisfying a Poincar\'e-type
inequality with constant $\sigma^2$, then, for all $h > 0$,
%
\begin{equation}\label{for4.4}
\mathbf{P}\biggl\{W_1(F_n,F) \geq
C \sigma\biggl(\frac{A + \log n}{n}\biggr)^{ 1/3} + h\biggr\}
\leq C \mathrm{e}^{- h\sqrt{n}/\sigma},
\end{equation}
where $A = \frac{1}{\sigma} \max_{i,j} | \mathbf{E}X_i - \mathbf{E}X_j|$
and where $C$ is an absolute constant.
\end{proposition}

Bounds such as (\ref{for4.4}) may be used to prove that the convergence holds
almost surely at a certain rate. Here is a simple example, corresponding
to non-varying values of the Poincar\'e constants. (One should
properly modify the conclusion when applying this to the matrix scheme;
see Section \ref{sec7}.) Let $(X_n)_{n \geq1}$ be a random sequence such that
for each $n$, $(X_1,\ldots,X_n)$ has distribution on $\mathbf{R}^n$ satisfying
$\operatorname{PI}(\sigma^2)$ with some common $\sigma$.

\begin{corollary}\label{cor4.3} If
$\max_{i,j \leq n} |\mathbf{E}X_i - \mathbf{E}X_j| = \mathrm{O}(\log n)$, then
$
W_1(F_n,F) = \mathrm{O}(\frac{\log n}{n})^{ 1/3}
$
with probability $1$.
\end{corollary}

Note, however, that in the scheme of sequences such as in Corollary \ref{cor4.3},
the mean distribution function $F = \mathbf{E}F_n$ might also depend
on $n$.

By a similar contraction argument, the upper bound (\ref{for4.4}) may be sharpened,
when the distribution of $(X_1,\ldots,X_n)$ satisfies a logarithmic
Sobolev inequality. We turn to this type of (stronger) hypothesis
in the next section.

\begin{remarks*} Let $\Omega$ be a metric space and let $d=d(u,u')$
be a non-negative continuous function on the product space
$\Omega\times\Omega$. Given Borel probability measures $G$ and $G'$
on $\Omega$, the generalized Kantorovich--Rubinstein or Wasserstein
`distance' with cost function $d$ is defined by
\[
W(G,G') = \inf_\pi\int d(u,u')\, \mathrm{d}\pi(u,u'),
\]
where the infimum is taken over all probability measures $\pi$
on $\Omega\times\Omega$ with marginal distributions~$G$ and $G'$.
In the case of the real line $\Omega= \mathbf{R}$ with cost function
of the
form $d(u,u') = \varphi(u - u')$, where $\varphi $ is convex, this quantity
has a simple description,
%
\begin{equation}\label{for4.5}
W(G,G') = \int_0^1 \varphi\bigl(G^{-1}(t) - G'^{-1}(t)\bigr)\, \mathrm{d}t,
\end{equation}
in terms of the inverse distribution functions
$G^{-1}(t) = \min\{x \in\mathbf{R}\dvtx G(x) \geq t\}$; see, for
example, \cite{C-S-S}
and \cite{R}, Theorem 2.

If $\varphi(u,u') = |u-u'|$, then we also have the $L^1$-representation
for $W_1(G,G')$, which we use from the very beginning as our definition.
Moreover, for arbitrary discrete measures
$G = F_n = (\delta_{x_1} + \cdots+ \delta_{x_n})/n$ and
$G' = F_n' = (\delta_{x_1'} + \cdots+ \delta_{x_n'})/n$,
as in Lemma \ref{lem4.1}, the expression (\ref{for4.5}) is reduced to
%
\begin{equation}\label{for4.6}
W_1(F_n,F_n') = \frac{1}{n} \sum_{i=1}^n |x_i - x_i'|,
\end{equation}
where we assume that $x_1 \leq\cdots\leq x_n$ and
$x_1' \leq\cdots\leq x_n'$.

Now, for an arbitrary random vector $X = (X_1,\ldots,X_n)$ in $\mathbf{R}^n$,
consider the
ordered statistics $X_1^* \leq\cdots\leq X_n^*$. Equation (\ref{for4.6}) then yields
%
\begin{equation}\label{for4.7}
\mathbf{E}W_1(F_n,F_n') = \frac{1}{n} \sum_{i=1}^n \mathbf{E}|X_i^*
- (X_i')^*|,
\end{equation}
where $(X_1')^* \leq\cdots\leq(X_n')^*$ are ordered statistics
generated by an independent copy of $X$ and where
$F_n'$ are independent copies of the (random) empirical measures $F_n$
associated with $X$. By the triangle inequality for the metric $W_1$,
we have
\[
\mathbf{E}W_1(F_n,F_n') \leq\mathbf{E}W_1(F_n,F) + \mathbf
{E}W_1(F,F_n') =
2 \mathbf{E}W_1(F_n,F).
\]
It is applied with the mean distribution function $F = \mathbf{E}F_n$.
On the other hand, any function of the form $H \rightarrow W_1(G,H)$
is convex on the convex set $Z_1$, so, by Jensen's inequality,
$\mathbf{E}W_1(F_n,F_n') \geq\mathbf{E}W_1(F_n,\mathbf{E}F_n') =
\mathbf{E}W_1(F_n,F)$.
The two bounds give
%
\begin{equation}\label{for4.8}
\mathbf{E}W_1(F_n,F) \leq\mathbf{E}W_1(F_n,F_n') \leq2 \mathbf
{E}W_1(F_n,F).
\end{equation}
By a similar argument,
%
\begin{equation}\label{for4.9}
\mathbf{E}|X_i^* - \mathbf{E}X_i^*| \leq\mathbf{E}|X_i^* -
(X_i')^*| \leq
2 \mathbf{E}|X_i^* - \mathbf{E}X_i^*|.
\end{equation}
Combining (\ref{for4.8}) and (\ref{for4.9}) and recalling (\ref{for4.7}), we arrive at the two-sided
estimate
\[
\frac{1}{2n} \sum_{i=1} \mathbf{E}|X_i^* - \mathbf{E}X_i^*| \leq
\mathbf{E}W_1(F_n,F)
\leq\frac{2}{n} \sum_{i=1} \mathbf{E}|X_i^* - \mathbf{E}X_i^*|,
\]
which is exactly the inequality (\ref{for1.4}) mentioned in the \hyperref[intr]{Introduction}.
Similar two-sided estimates also hold for other cost
functions in the Wasserstein distance.
\end{remarks*}

\section{Empirical log-Sobolev inequalities}\label{sec5}

As before, let $(X_1,\ldots,X_n)$ be a random vector in $\mathbf
{R}^n$ with joint
distribution $\mu$. Similarly to Proposition \ref{pro2.1}, now using a
log-Sobolev inequality for $\mu$, we arrive at the following `empirical'
log-Sobolev inequality.

\begin{proposition}\label{pro5.1} Under $\operatorname{LSI}(\sigma^2)$,
for any bounded, smooth function $f$ on $\mathbf{R}$,
\[
\operatorname{Ent}_\mu\biggl[\biggl(\int f\,\mathrm{d}F_n\biggr)^{ 2}\biggr] \leq
\frac{2\sigma^2}{n} \int f'^2\,\mathrm{d}F.
\]
\end{proposition}

In analogy with Poincar\'e-type inequalities, one may also develop refined
applications to the rate of growth of moments and to large deviations
of various functionals of empirical measures.
In particular, we have the following proposition.

\begin{proposition}\label{pro5.2} Under $\operatorname{LSI}(\sigma^2)$, for
any smooth function $f$ on $\mathbf{R}$ such that $f'$ belongs to
$L^p(\mathbf{R},\mathrm{d}F)$, $p \geq2$,
%
\begin{equation}\label{for5.1}
\mathbf{E}\biggl|\int f\,\mathrm{d}F_n - \int f\,\mathrm{d}F \biggr|^p \leq
\frac{(\sigma\sqrt{p} )^p}{n^{p/2}} \int|f'|^p\,\mathrm{d}F.
\end{equation}
In addition, if $|f'| \leq1$, then, for all $h>0$,
%
\begin{equation}\label{for5.2}
\mu\biggl\{\biggl|\int f\,\mathrm{d}F_n - \int f\,\mathrm{d}F \biggr| \geq h\biggr\}
\leq2 \mathrm{e}^{-nh^2/2\sigma^2}.
\end{equation}
\end{proposition}

The proof of the second bound, (\ref{for5.2}), which was already noticed in
\cite{G-Z} in the context of random matrices, follows the standard Herbst's
argument; see \cite{L1} and \cite{B-G1}.
The first family of moment inequalities, (\ref{for5.1}), can be sharpened
by one inequality on the Laplace transform, such as
\[
\mathbf{E}\exp\biggl\{\int f\,\mathrm{d}F_n - \int f\,\mathrm{d}F\biggr\} \leq
\mathbf{E}\exp\biggl\{\frac{\sigma^2}{n} \int|f'|^2\,\mathrm{d}F_n\biggr\}.
\]
The proof is immediate, by \cite{B-G1}, Theorem 1.2.

However, a major weak point in both Poincar\'e and log-Sobolev inequalities,
including their direct consequences, as in Proposition \ref{pro5.2}, is that
they may not be applied to indicator and other non-smooth functions.
In particular, we cannot estimate directly at fixed points the
variance $\operatorname{Var}(F_n(x))$ or other similar quantities like
the higher
moments of $|F_n(x) - F(x)|$. Therefore, we need another family of
analytic inequalities. Fortunately, the so-called infimum-convolution
operator and associated relations concerning arbitrary measurable functions
perfectly fit our purposes. Moreover, some of the important relations
hold true and may be controlled in terms of the constant involved
in the logarithmic Sobolev inequalities.

Let us now turn to the important concept of infimum- and
supremum-convolution inequalities. They were proposed in 1991 by
Maurey \cite{Mau} as a functional approach to some of Talagrand's
concentration results concerning product measures.
Given a parameter $t > 0$ and a real-valued function $g$ on $\mathbf{R}^n$
(possibly taking the values $\pm\infty$), put
\begin{eqnarray*}
Q_t g(x) &=& \inf_{y \in\mathbf{R}^n} \biggl[ g(y) + \frac{1}{2t} \|x-y\|
^2\biggr],
\\
P_t g(x) &=& \sup_{y \in\mathbf{R}^n} \biggl[ g(y) - \frac{1}{2t} \|x-y\|
^2\biggr].
\end{eqnarray*}
$Q_t g$ and $P_t g$ then represent, respectively, the infimum- and
supremum-convolution of $g$ with cost function being the
normalized square of the Euclidean norm in $\mathbf{R}^n$.
By definition, one puts $Q_0 g = P_0 g = g$.

For basic definitions and basic properties of the infimum- and
supremum-convolution operators, we refer the reader to \cite{Ev} and
\cite{B-G-L}, mentioning just some of them here. These operators are dually
related by
the property that for any functions
$f$ and $g$ on $\mathbf{R}^n$, $g \geq P_t f \Longleftrightarrow f
\leq Q_t g$.
Clearly, $P_t (-g) = -Q_t g$. Thus, in many statements, it is
sufficient to consider only one of these operators.
The basic semigroup property of both operators is that
for any $g$ on $\mathbf{R}^n$ and $t,s \geq0$,
\[
Q_{t+s} g = Q_t Q_s g,\qquad  P_{t+s} g = P_t P_s g.
\]

For any function $g$ and $t>0$, the function $P_t g$ is always lower
semicontinuous, while $Q_t g$ is upper semicontinuous. If $g$ is bounded,
then $P_t g$ and $Q_t g$ are bounded and have finite Lipschitz seminorms.
In particular, both are differentiable almost everywhere.

Given a bounded function $g$ and $t>0$, for almost all $x \in\mathbf{R}^n$,
the functions $t \rightarrow P_t g(x)$ and $t \rightarrow Q_t g(x)$
are differentiable at $t$ and
\[
\frac{\partial P_t g(x)}{\partial t} = \frac{1}{2} \|\nabla P_t g(x)\|
^2,
\qquad\frac{\partial Q_t g(x)}{\partial t} = - \frac{1}{2} \|\nabla Q_t
g(x)\|^2\qquad
\mbox{a.e.}
\]
In other words, the operator $\Gamma g = \frac{1}{2} |\nabla g|^2$
appears as the generator for the semigroup $P_t$, while~${-}\Gamma$
appears as the generator for $Q_t$. As a result, $u(x,t) = Q_t g(x)$ represents
the solution to the Hamilton--Jacobi equation
$\frac{\partial u}{\partial t} = - \frac{1}{2} \|\nabla u\|^2$
with initial condition $u(x,0) = g(x)$.

Below, we separately formulate a principal result of \cite{B-G-L} which
relates logarithmic Sobolev inequalities to supremum- and
infimum-convolution operators.

\begin{lemma}\label{lem5.3} Let $\mu$ be a probability measure on
$\mathbf{R}^n$
satisfying $\operatorname{LSI}(\sigma^2)$. For any
$\mu$-integrable Borel-measurable function $g$ on $\mathbf{R}^n$,
we have
%
\begin{equation}\label{for5.3}
\int P_{\sigma^2} g \,\mathrm{d}\mu\geq\log\int \mathrm{e}^g \,\mathrm{d}\mu
\end{equation}
and, equivalently,
%
\begin{equation}\label{for5.4}
\int g \,\mathrm{d}\mu\geq\log\int \mathrm{e}^{Q_{\sigma^2} g} \,\mathrm{d}\mu.
\end{equation}
\end{lemma}

Alternatively, for further applications to empirical measures, one could
start from the infimum-convolution inequalities (\ref{for5.3}) and (\ref{for5.4}), taking them
as the main hypothesis on the measure $\mu$. They take an intermediate
position between Poincar\'e and logarithmic Sobolev inequalities.
However, logarithmic Sobolev inequalities have been much better studied,
with a variety of sufficient conditions having been derived.

Now, as in Section \ref{sec2}, we apply the relations (\ref{for5.3}) and (\ref{for5.4}) to functions
$g(x_1,\ldots,x_n) = \int f\,\mathrm{d}F_n$,
where $F_n$ is the empirical measure defined for `observations'
$x_1,\ldots,x_n$. By the very definition, for any $t>0$,
\begin{eqnarray*}
P_t g(x_1,\ldots,x_n) & = &
\sup_{y_1,\ldots,y_n \in\mathbf{R}} \Biggl[ g(y_1,\ldots,y_n) -
\frac{1}{2t} \sum_{i=1}^n |x_i-y_i|^2\Biggr] \\
& = &
\frac{1}{n} \sup_{y_1,\ldots,y_n \in\mathbf{R}} \sum_{i=1}^n
\biggl[ f(y_i) - \frac{1}{2t/n} |x_i-y_i|^2\biggr]  =  \int P_{t/n} f\,\mathrm{d}F_n.
\end{eqnarray*}
Similarly, $Q_t g = \int Q_{t/n} f\,\mathrm{d}F_n$. Therefore, after integration
with respect to $\mu$ and using the identity $P_a (tf) = t P_{ta} f$,
we arrive at corresponding empirical supremum- and infimum-convolution
inequalities, as follows.

\begin{proposition}\label{pro5.4} Under $\operatorname{LSI}(\sigma^2)$, for any
$F$-integrable Borel-measurable function $f$ on $\mathbf{R}$ and for
any $t>0$,
%
\begin{eqnarray}\label{for5.5}
\log\mathbf{E}\mathrm{e}^{t (\int f\,\mathrm{d}F_n - \int f\,\mathrm{d}F)} & \leq&
t \int[P_{t\sigma^2/n} f - f ]\,\mathrm{d}F,  \\
\label{for5.6}\log\mathbf{E}\mathrm{e}^{t (\int f\,\mathrm{d}F - \int f\,\mathrm{d}F_n)} & \leq&
t \int[f - Q_{t\sigma^2/n} f ]\,\mathrm{d}F.
\end{eqnarray}
\end{proposition}

Note that the second inequality may be derived from the first by changing
$f$ to $-f$.

\section{Local behavior of empirical distributions}\label{sec6}

In this section, we develop a few direct applications of Proposition \ref{pro5.4}
to the behavior of empirical distribution functions $F_n(x)$ at a fixed
point. Such functionals are linear, that is, of the form $\int f\,\mathrm{d}F_n$,
corresponding to the indicator function of the half-axis $f = (-\infty,x]$.
When $f$ is smooth, Proposition \ref{pro5.2} tells us that the deviations of
$L_n f = \int f\,\mathrm{d}F_n - \int f\,\mathrm{d}F$ are of order $\sigma/\sqrt{n}$.
In the general non-smooth case, the infimum- and supremum-convolution
operators $P_t f$ and $Q_t f$ behave differently for small values of
$t$ and this results in a different rate of fluctuation for $L_n f$.

To see this, let $f = (-\infty,x]$. In this case, the functions
$P_t f$ and $Q_t f$ may easily be computed explicitly, but we do not
lose much by using the obvious bounds
\[
1_{(-\infty, x - \sqrt{2t} ]} \leq Q_t f \leq
P_t f \leq1_{(-\infty, x + \sqrt{2t} ]}.
\]
Therefore, (\ref{for5.5}) and (\ref{for5.6}) yield the following proposition.

\begin{proposition}\label{pro6.1} Under $\operatorname{LSI}(\sigma^2)$, for any
$x \in\mathbf{R}$ and $t>0$, with $h = \sqrt{\frac{2\sigma^2 t}{n}}$,
%
\begin{eqnarray}\label{for6.1}
\log\mathbf{E}\mathrm{e}^{t (F_n(x) - F(x))} & \leq&
t \bigl(F(x + h) - F(x)\bigr),  \\
\label{for6.2}\log\mathbf{E}\mathrm{e}^{t (F(x) - F_n(x))} & \leq&
t \bigl(F(x) - F(x-h)\bigr).
\end{eqnarray}
\end{proposition}

These estimates may be used to sharpen Corollary \ref{cor3.2} and therefore
to recover Theorem~\ref{teo1.1} (under the stronger hypothesis on
the joint distribution $\mu$, however). Indeed, for any $t>0$,
\[
\mathbf{E}\mathrm{e}^{t |F_n(x) - F(x)|} \leq
\mathbf{E}\mathrm{e}^{t (F_n(x) - F(x))} + \mathbf{E}\mathrm{e}^{-t (F(x) - F_n(x))}
\leq
2 \mathrm{e}^{t (F(x+h) - F(x-h))}.
\]
Taking the logarithm and applying Jensen's inequality, we arrive at
\[
\mathbf{E}|F_n(x) - F(x)| \leq\bigl(F(x+h) - F(x-h)\bigr) + \frac{\log2}{t}.
\]
Now, just integrate this inequality over an arbitrary interval $(a,b)$,
$a<b$, and use the general relation
$\int_{-\infty}^{+\infty} (F(x+h) - F(x-h)) \,\mathrm{d}x \leq2h$
to obtain that
\[
\mathbf{E}\int_a^b |F_n(x)-F(x)| \,\mathrm{d}x \leq2h + \frac{\log2}{t} (b-a)
= 2\sqrt{\frac{2\sigma^2 t}{n}} + \frac{\log2}{t} (b-a).
\]
Optimization over $t$ leads to an improved version of Corollary \ref{cor3.2}.

\begin{corollary}\label{cor6.2} Under $\operatorname{LSI}(\sigma^2)$, for all $a<b$,
\[
\mathbf{E}\int_a^b |F_n(x)-F(x)| \,\mathrm{d}x \leq4
\biggl(\frac{\sigma^2 (b-a)}{n}\biggr)^{ 1/3}.
\]
\end{corollary}

Note that in both cases of Proposition \ref{pro6.1}, for any $t \in\mathbf{R}$,
\[
\log\mathbf{E}\mathrm{e}^{t (F_n(x) - F(x))} \leq|t| \bigl(F(x+h) - F(x-h)\bigr),\qquad
h = \sqrt{2\sigma^2 |t|/n}.
\]
Hence, the local behavior of the distribution function $F$ near
a given point $x$ turns out to be responsible for the large deviation
behavior at this point of the empirical distribution function $F_n$
around its mean.

For a quantitative statement, assume that $F$ has a finite Lipschitz
constant $M = \|F\|_{\mathrm{Lip}}$, so it is absolutely continuous
with respect to Lebesgue measure on the real line and has a density,
bounded by $M$. It follows from (\ref{for6.1}), with
$t = (\alpha n^{1/3}) \lambda$ and $\alpha^3 = \frac{2}{9 M^2 \sigma^2}$,
that
\[
\mathbf{E}\mathrm{e}^{\lambda\xi} \leq \mathrm{e}^{({2}/{3}) |\lambda|^{3/2}},
\qquad\lambda\in\mathbf{R},
\]
where $\xi= \alpha n^{1/3} (F_n(x) - F(x))$. By Chebyshev's
inequality, for any $r > 0$,
\[
\mu\{\xi\geq r\} \leq \mathrm{e}^{({2}/{3}) \lambda^{3/2} - \lambda r} =
\mathrm{e}^{-r^3/3},\qquad   \mbox{where } \lambda= r^2.
\]
Similarly, $\mu\{\xi\leq-r\} \leq \mathrm{e}^{-r^3/3}$. Therefore,
$\mu\{\alpha n^{1/3} |F_n(x) - F(x)| \geq r\} \leq2 \mathrm{e}^{-r^3/3}$.
Changing the variable, we are finished.

Recall that we use the quantity
\[
\beta= \frac{(M\sigma)^{2/3}}{n^{1/3}}.
\]

\begin{proposition}\label{pro6.3} Assume that $F$ has a density,
bounded by
a number $M$. Under $\operatorname{LSI}(\sigma^2)$, for any $x \in\mathbf{R}$
and $r>0$,
%
\begin{equation}\label{for6.3}
\mathbf{P}\{|F_n(x) - F(x)| \geq\beta r\} \leq2 \mathrm{e}^{-2r^3/27}.
\end{equation}
In particular, with some absolute constant $C$, we have
%
\begin{equation}\label{for6.4}
\mathbf{E}|F_n(x) - F(x)| \leq C \beta.
\end{equation}
\end{proposition}

Note that (\ref{for6.4}) is consistent with the estimate of Theorem~\ref{teo1.1}.
To derive similar bounds on the uniform (Kolmogorov) distance
$\|F_n-F\| = \sup_x |F_n(x)-F(x)|$ (which we discuss in the next section),
it is better to split the bound (\ref{for6.3}) into the two parts,
%
\begin{eqnarray}\label{for6.5}
\mathbf{P}\{F_n(x) - F(x) \geq\beta r\} & \leq& \mathrm{e}^{-2r^3/27},  \\
\label{for6.6}\mathbf{P}\{F(x) - F_n(x) \geq\beta r\} & \leq& \mathrm{e}^{-2r^3/27},
\end{eqnarray}
which were obtained in the last step of the proof of (\ref{for6.3}).

However, since one might not know whether $F$ is Lipschitz or how it
behaves locally, and since one might want to approximate this measure
itself by some canonical distribution $G$, it is reasonable to provide
a more general statement. By Proposition \ref{pro6.1}, for any $t \geq0$,
\begin{eqnarray*}
\log\mathbf{E}\mathrm{e}^{t (F_n(x) - G(x))} & \leq& t \bigl(F(x+h) - G(x)\bigr) \\
& \leq&
t \bigl(G(x + h) - G(x)\bigr) + t \|F-G\|
\end{eqnarray*}
and, similarly,
\[
\log\mathbf{E}\mathrm{e}^{t (G(x) - F_n(x))} \leq
t \bigl(G(x) - G(x-h)\bigr) + t \|F-G\|.
\]
Repeating the preceding argument with the random variable
\[
\xi= \alpha n^{1/3} \bigl(F_n(x) - G(x) - \|F-G\|\bigr)
\]
and then interchanging $F_n$ and $G$, we get a more general version
of Proposition \ref{pro6.3}.

\begin{proposition}\label{pro6.4} Under $\operatorname{LSI}(\sigma^2)$, for any
distribution function $G$ with finite Lipschitz seminorm
$M = \|G\|_{\mathrm{Lip}}$, for any $x \in\mathbf{R}$ and $r>0$,
\[
\mathbf{P}\{|F_n(x) - G(x)| \geq\beta r + \|F-G\| \} \leq
2 \mathrm{e}^{-2r^3/27},
\]
where $\beta= (M\sigma)^{2/3} n^{-1/3}$.
In particular, up to some absolute constant $C$,
%
\begin{equation}\label{for6.7}
\mathbf{E}|F_n(x) - G(x)| \leq C \beta+ \|F-G\|.
\end{equation}
\end{proposition}

Let us stress that in all of these applications of Proposition \ref{pro5.4}, only
the indicator functions $f = 1_{(-\infty,x]}$ were used. One may
therefore try
to get more information about deviations of the empirical distributions
$F_n$ from the mean $F$ by applying the basic bounds (\ref{for5.5}) and (\ref{for5.6})
with different (non-smooth) functions $f$.

For example, of considerable interest is the so-called local regime, where
one tries to estimate the number
\[
N_I = \operatorname{card}\{i \leq n\dvtx X_i \in I\}
\]
of observations inside a small interval $I = [x,x+\varepsilon ]$ and
to take into account the size of the increment $\varepsilon = |I|$.
In case of i.i.d. observations, this may done using various tools;
already, the formula
\[
\operatorname{Var}(F_n(I)) = \frac{1}{n} F(I) \bigl(1 - F(I)\bigr) \leq
\frac{F(x+\varepsilon ) - F(x)}{n}
\]
suggests that when $F$ is Lipschitz, $F_n(I)$ has small oscillations
for small $\varepsilon $ (where $F_n$ and $F$ are treated as measures).

However, the infimum- and supremum-convolution operators $P_t f$ and
$Q_t f$ do not provide such information. Indeed, for the indicator
function $f = 1_I$, by (\ref{for5.5}), we only have, similarly to
Proposition \ref{pro6.1}, that
\[
\log\mathbf{E}\mathrm{e}^{t (F_n(I) - F(I))} \leq
t \bigl[\bigl(F(x + \varepsilon + h) - F(x)\bigr) - \bigl(F(x) - F(x-h)\bigr)\bigr],
\]
where $t>0$ and $h = \sqrt{\frac{2\sigma^2 t}{n}}$. Here, when $h$ is
fixed and $\varepsilon \rightarrow0$, the right-hand side is not vanishing,
in contrast with the i.i.d. case. This also shows that standard chaining
arguments, such as Dudley's entropy bound or more delicate majorizing
measure techniques (described, e.g., in \cite{Ta}),
do not properly work through the infimum-convolution approach.

Nevertheless, the above estimate is still effective for $\varepsilon $
of order $h$,
so we can control deviations of $F_n(I) - F(I)$ relative to
$|I|$ when the intervals are not too small. This can be done with the
arguments used in the proof of Proposition \ref{pro6.3} or, alternatively
(although with worse absolute constants), one can use the inequality (\ref{for6.3}),
by applying it to the points $x$ and $x+\varepsilon $. This
immediately gives that
\[
\mathbf{P}\{|F_n(I) - F(I)| \geq2\beta r\} \leq
4 \mathrm{e}^{-2r^3/27}.
\]
Changing variables, one may rewrite the above in terms of $N_I$ as
\[
\mathbf{P}\{|N_I - n F(I) | \geq n\delta|I|\} \leq
4 \exp\biggl\{-c \biggl(\frac{\delta|I|}{\beta}\biggr)^3\biggr\}
\]
with $c = 1/112$. Note that the right-hand side is small only when
$|I| \gg\beta/\delta$, which is of order $n^{-1/3}$ with respect to
the number of observations.

This can further be generalized if we apply Proposition \ref{pro6.4}.

\begin{corollary}\label{cor6.5} Let $G$ be a distribution function with
density $g(x)$ bounded by a number $M$. Under $\operatorname{LSI}(\sigma^2)$,
for any
$\delta> 0$ and any interval $I$ of length $|I| \geq4 \|F-G\|/\delta$,
\[
\mathbf{P}\biggl\{\biggl|N_I - n \int_I g(x) \,\mathrm{d}x \biggr|
\geq n\delta|I|\biggr\} \leq4
\exp\biggl\{ -c \biggl(\frac{\delta|I|}{\beta}\biggr)^{ 3}\biggr\},
\]
where $\beta= (M\sigma)^{2/3} n^{-1/3}$ and $c>0$ is an absolute
constant.
\end{corollary}

Hence, if
\[
|I| \geq\frac{C}{\delta} \max\{\beta,\|F-G\|\}
\]
and $C>0$ is large, then, with high probability, we have that
$|\frac{N_I}{n} - \int_I g(x) \,\mathrm{d}x | \leq\delta|I|$.

\section[Bounds on the Kolmogorov distance. Proof of Theorem 1.2]{Bounds on
the Kolmogorov distance. Proof of Theorem \protect\ref{teo1.2}}\label{sec7}

As before, let $F_n$ denote the empirical measure associated with
observations $x_1,\ldots,x_n$ and $F = \mathbf{E}F_n$ their mean with respect
to a given probability measure $\mu$ on $\mathbf{R}^n$. In this
section, we
derive uniform bounds on $F_n(x)-F(x)$, based on Proposition \ref{pro6.3},
and thus prove Theorem \ref{teo1.2}.
For applications to the matrix scheme, we shall also replace $F$,
which may be difficult to determine, by the well-behaving limit law $G$
(with the argument relying on Proposition \ref{pro6.4}).

Let the random variables $X_1,\ldots,X_n$ have joint distribution $\mu$,
satisfying $\operatorname{LSI}(\sigma^2)$, and assume that $F$ has a finite Lipschitz
seminorm $M = \|F\|_{\mathrm{Lip}}$. Define
\[
\beta= \frac{(M\sigma)^{2/3}}{n^{1/3}}.
\]

\begin{pf*}{Proof of Theorem \ref{teo1.2}}
We use the inequalities (\ref{for6.5}) and (\ref{for6.6}) to derive an upper bound on
$\|F_n - F\| = \sup_x |F_n(x)-F(x)|$.
(For the sake of extension of Theorem \ref{teo1.2} to Theorem \ref{teo7.1} below, we relax
the argument and do not assume that $F$ is continuous.)

So, fix $r > 0$ and an integer
$N \geq2$. One can always pick up points
$-\infty= x_0 \leq x_1 \leq\cdots\leq x_{N-1} \leq x_N = +\infty$
with the property that
%
\begin{equation}\label{for7.1}
F(x_i-) - F(x_{i-1}) \leq\frac{1}{N}, \qquad i = 1,\ldots,N.
\end{equation}
Note that $F_n(x_0) = F_n(x_0) = 0$ and $F_n(x_N) = F_n(x_N) = 1$.
It then follows from (\ref{for6.5}) that
\[
\mathbf{P}\Bigl\{\max_{1 \leq i \leq N} [F_n(x_i-) - F(x_i-)] \geq\beta
r\Bigr\}
\leq(N-1) \mathrm{e}^{-2r^3/27}
\]
and, similarly, by (\ref{for6.6}),
\[
\mathbf{P}\Bigl\{\max_{1 \leq i \leq N} [F(x_{i-1}) - F_n(x_{i-1})] \geq
\beta r\Bigr\}
\leq(N-1) \mathrm{e}^{-2r^3/27}.
\]
Hence, for the random variable
\[
\xi_N = \max\Bigl\{
\max_{1 \leq i \leq N} [F_n(x_i-) - F(x_i-)],
\max_{1 \leq i \leq N} [F(x_{i-1}) - F_n(x_{i-1})]\Bigr\},
\]
we have that
%
\begin{equation}\label{for7.2}
\mathbf{P}\{\xi_N \geq\beta r\} \leq2(N-1) \mathrm{e}^{-2r^3/27}.
\end{equation}

Now, take any point $x \in\mathbf{R}$ different from all of the
$x_j$'s and select
$i$ from $1,\ldots,n$ such that $x_{i-1} < x < x_i$. Then, by (\ref{for7.1}),
\begin{eqnarray*}
F_n(x) - F(x) & \leq& F_n(x_i-) - F(x_{i-1}) \\
& = & [F_n(x_i-) - F(x_i-)] + [F(x_i-) - F(x_{i-1})]
 \leq \xi_N + \frac{1}{N}.
\end{eqnarray*}
Similarly,
\begin{eqnarray*}
F(x) - F_n(x) & \leq& F(x_i-) - F_n(x_{i-1}) \\
& = & [F(x_i-) - F(x_{i-1})] + [F(x_{i-1}) - F_n(x_{i-1})]
 \leq \xi_N + \frac{1}{N}.
\end{eqnarray*}
Therefore, $|F_n(x) - F(x)| \leq\xi_N + \frac{1}{N}$, which also
extends by continuity from the right to all points $x_j$.
Thus, $\|F_n - F\| \leq\xi_N + \frac{1}{N}$ and, by (\ref{for7.2}),
\[
\mathbf{P}\biggl\{\|F_n - F\| > \beta r + \frac{1}{N}\biggr\} \leq
2(N-1) \mathrm{e}^{-2r^3/27}.
\]
Note that this also holds automatically in the case $N=1$.
Choose $N = [\frac{1}{\beta r}] + 1$. We then have $\frac{1}{N} \leq
\beta r$ and get
\[
\mathbf{P}\{\|F_n - F\| > 2\beta r \} \leq
\frac{2}{\beta r} \mathrm{e}^{-2r^3/27}.
\]
Finally, changing $2\beta r$ into $r$, we arrive at the bound (\ref{for1.5}) of
Theorem \ref{teo1.2},
%
\begin{equation}\label{for7.3}
\mathbf{P}\{\|F_n - F\| > r \} \leq
\frac{4}{r} \exp\biggl\{-\frac{2}{27} \biggl(\frac{r}{\beta}\biggr)^3\biggr\},
\end{equation}
and so the constant $c = 2/27$.

It remains to derive the bound on the mean $\mathbf{E}\|F_n - F\|$. Given
$0 \leq r_0 \leq1$, we can write, using~(\ref{for7.3}),
%
\begin{eqnarray}\label{for7.4}
\mathbf{E}\|F_n - F\| & = & \int_0^1 \mu\{\|F_n - F\| > r\} \,\mathrm{d}r  =  \int_0^{r_0} + \int_{r_0}^1
\nonumber
\\[-8pt]
\\[-8pt]
\nonumber
& \leq& r_0 + \frac{4}{r_0} \exp\biggl\{-\frac{2}{27} \biggl(\frac{r_0}{\beta
}\biggr)^3\biggr\}.
\end{eqnarray}
This bound also holds for $r_0>1$.

First, assume that $0 < \beta\leq1$ and choose
$r_0 = 3\beta\log(1 + \frac{1}{\beta})$. Then, for the last term
in (\ref{for7.4}), we have
\begin{eqnarray*}
\frac{4}{r_0} \exp\biggl\{-\frac{2}{27} \biggl(\frac{r_0}{\beta}\biggr)^3\biggr\}
& = &
\frac{4}{3\beta\log(1 + {1}/{\beta})}  \mathrm{e}^{-2 \log(1 + {1}/{\beta})} \\
& = &
\frac{4}{3} \frac{\beta}{(1+\beta)^2 \log^{1/3}(1 + {1}/{\beta})}\\
 &\leq& B\beta\log^{1/3}\biggl(1 + \frac{1}{\beta}\biggr)
\end{eqnarray*}
with some constant $B$ satisfying
$(1+\beta)^3 \log(1 + \frac{1}{\beta}) \geq(\frac{4}{3 B})^{3/2}$.
For example, we can take $B=2$ and then, by (\ref{for7.4}), we have
%
\begin{equation}\label{for7.5}
\mathbf{E}\|F_n - F\| \leq5 \beta\log^{1/3}\biggl(1 + \frac{1}{\beta}\biggr).
\end{equation}

As for the values $\beta\geq1$, simple calculations show that
the right-hand side of (\ref{for7.5}) is greater than 1, so the inequality
(\ref{for1.6}) is fulfilled with $C = 5$.

Theorem \ref{teo1.2} is therefore proved.
\end{pf*}

\begin{remark*} If $M\sigma$ in Theorem 1.2 were of order 1, then
$\mathbf{E}\|F_n - F\|$ would be of order at most
$(\frac{\log n}{n})^{1/3}$.
Note, however, that under $\operatorname{PI}(\sigma^2)$, and if all $\mathbf
{E}X_i = 0$,
the quantity $M\sigma$ is separated from zero and, more precisely,
$M \sigma\geq\frac{1}{\sqrt{12}}$.
\end{remark*}

Indeed, by Hensley's theorem in dimension $1$ \cite{H,Ba}, in the class
of all probability densities~$p(x)$ on the line, the expression
$(\int x^2 p(x) \,\mathrm{d}x)^{1/2} \operatorname{ess\, sup}_x p(x)$ is minimized for
the uniform distribution on symmetric intervals and is therefore bounded
from below by $1/\sqrt{12}$. Since $F$ is Lipschitz, it has a density
$p$ with $M = \operatorname{ess\, sup}_x p(x)$. On the other hand, it follows from
the Poincar\'e-type inequality that $\sigma^2 \geq\operatorname
{Var}(X_i) = \mathbf{E}X_i^2$.
Averaging over all $i$'s, we get $\sigma^2 \geq\int x^2\,\mathrm{d}F(x)$, so
$M\sigma\geq(\int x^2 p(x) \,\mathrm{d}x)^{1/2} \operatorname{ess\, sup}_x p(x)$.

With similar arguments based on Proposition \ref{pro6.4}, we also obtain the
following generalization of Theorem \ref{teo1.2}.

\begin{theorem}\label{teo7.1} Assume that $X_1,\ldots,X_n$ have a
distribution
on $\mathbf{R}^n$ satisfying $\operatorname{LSI}(\sigma^2)$. Let $G$ be
a distribution function with finite Lipschitz seminorm $M$. Then,
for all $r > 0$,
\[
\mathbf{P}\{\|F_n - G\| \geq r + \|F - G\|\} \leq
\frac{4}{r} \exp\biggl\{-\frac{2}{27} \biggl(\frac{r}{\beta}\biggr)^3\biggr\},
\]
where $\beta= (M\sigma)^{2/3} n^{-1/3}$. In particular,
\[
\mathbf{E}\|F_n - G\| \leq5 \beta\log^{1/3}\biggl(1 + \frac{1}{\beta}\biggr)
+ \|F-G\|.
\]
\end{theorem}

\begin{remarks*} It remains unclear whether or not one can involve
the i.i.d. case in the scheme of Poincar\'e or logarithmic Sobolev
inequalities to recover the rate $1/\sqrt{n}$ for
$\mathbf{E}\|F_n - F\|$, even if some further natural assumptions are
imposed (which are necessary, as we know from Examples~\ref{ex1} and~\ref{ex2}).
In particular, one may assume that the quantities $M$ and $\sigma$
are of order $1$, and that $\mathbf{E}X_i = 0$, $\operatorname
{Var}(X_i) = 1$. The question
is the following:
under, say, $\operatorname{LSI}(\sigma^2)$, is it true that
\[
\mathbf{E}\|F_n - F\| \leq\frac{C}{\sqrt{n}}
\]
with some absolute $C$? Or at least
$\mathbf{E}|F_n(x) - F(x)| \leq\frac{C}{\sqrt{n}}$ for individual points?
\end{remarks*}


\section{High-dimensional random matrices}\label{sec8}
\setcounter{equation}{0}

We shall now apply the bounds obtained in Theorems \ref{teo1.1} and \ref{teo7.1}, to the
case of the spectral empirical distributions. Let
$\{\xi_{jk}\}_{1 \leq j \leq k \leq n}$ be a family of independent
random variables on some probability space with mean $\mathbf{E}\xi
_{jk} = 0$
and variance $\operatorname{Var}(\xi_{jk}) = 1$. Put $\xi_{jk} = \xi
_{kj}$ for
$1 \leq k < j \leq n$ and introduce a symmetric $n \times n$ random
matrix,
\[
\Xi =  \frac{1}{\sqrt{n}}
\pmatrix{
\xi_{11} & \xi_{12} & \cdots& \xi_{1n} \cr
\xi_{21} & \xi_{22} & \cdots& \xi_{2n} \cr
\vdots& \vdots& \ddots& \vdots\cr
\xi_{n1} & \xi_{n2} & \cdots& \xi_{nn}}.
\]
Arrange its (real random) eigenvalues in increasing order:
$X_1 \leq\cdots\leq X_n$. As before, we associate with particular values
$X_1 = x_1, \ldots, X_n = x_n$ an empirical (spectral) measure $F_n$
with mean (expected) measure $F = \mathbf{E}F_n$.

An important point in this scheme is that the joint distribution $\mu$
of the spectral values, as a probability measure on $\mathbf{R}^n$,
represents the image of the joint distribution of $\xi_{jk}$'s
under a Lipschitz map $T$ with Lipschitz seminorm
$\|T\|_{\mathrm{Lip}} = \frac{\sqrt{2}}{\sqrt{n}}$. More precisely,
by the Hoffman--Wielandt theorem with respect to the Hilbert--Schmidt
norm, we have
\[
\sum_{i=1}^n |X_i - X_i'|^2 \leq\|\Xi- \Xi'\|_{\mathrm{HS}}^2
= \frac{1}{n} \sum_{j,k=1}^n |\xi_{jk} - \xi_{jk}'|^2 \leq
\frac{2}{n} \sum_{1 \leq j \leq k \leq n} |\xi_{jk} - \xi_{jk}'|^2
\]
for any collections $\{\xi_{jk}\}_{j \leq k}$ and
$\{\xi_{jk}\}_{j \leq k}'$ with eigenvalues $(X_1,\ldots,X_n)$,
$(X_1',\ldots,X_n')$, respectively.
This is a well-known fact (\cite{Bh}, page 165) which may be used
in concentration problems; see, for example, \cite{L2,D-S}.

In particular (see Proposition \ref{pra1} in the \hyperref[app]{Appendix}), if the
distributions of $\xi_{jk}$'s satisfy a one-dimensional Poincar\'e-type
inequality with common constant $\sigma^2$, then $\mu$ satisfies
a Poincar\'e-type inequality with an asymptotically much better constant
$\sigma_n^2 = \frac{2\sigma^2}{n}$. According to Theorem~\ref{teo1.1},
\[
\mathbf{E}\int_{-\infty}^{+\infty} |F_n(x) - F(x)| \,\mathrm{d}x \leq
C\sigma_n \biggl(\frac{A_n + \log n}{n}\biggr)^{ 1/3},
\]
where $C$ is an absolute constant and
$A_n = \frac{1}{\sigma_n} \max_{i,j} |\mathbf{E}X_i - \mathbf{E}X_j|$.
Since $\max_i |\mathbf{E}X_i|$ is of order at most $\sigma$, $A_n$
is at most
$\sqrt{n}$ and we arrive at the bound (\ref{for1.7}) in Theorem \ref{teo1.3}:
\[
\mathbf{E}\int_{-\infty}^{+\infty} |F_n(x) - F(x)| \,\mathrm{d}x \leq
\frac{C\sigma}{n^{2/3}}.
\]

Now, let us explain the second statement of Theorem \ref{teo1.3} for the case
where the $\xi_{jk}$'s satisfy a logarithmic Sobolev inequality
with a common constant $\sigma^2$, in addition to the normalizing
conditions $\mathbf{E}\xi_{jk} = 0$, $\operatorname{Var}(\xi_{jk})
= 1$ (which implies that
$\sigma\geq1$). Let $G$ denote the
standard semicircle law with variance 1, that is, with density
$g(x) = \frac{1}{2\curpi} \sqrt{4 - x^2}$, $-2 < x < 2$. In this case,
the Lipschitz seminorm is $M = \|G\|_{\mathrm{Lip}} = \frac{1}{\curpi}$. Also,
\[
\beta_n = \frac{(M\sigma_n)^{2/3}}{n^{1/3}} =
C' \biggl(\frac{\sigma}{n}\biggr)^{ 2/3}
\]
for some absolute $C'$. Therefore, applying Theorem \ref{teo7.1}
and using $\sigma\geq1$, we arrive at the bound~(\ref{for1.8}):
%
\begin{equation}\label{for8.1}
\mathbf{E}\sup_{x \in\mathbf{R}^n} |F_n(x) - G(x)| \leq C\sigma^{2/3}
\frac{\log^{1/3} n}{n^{2/3}} + \sup_{x \in\mathbf{R}} |F(x)-G(x)|.
\end{equation}

Thus, Theorem \ref{teo1.3} is proved.
For individual points that are close to the end-points $x = \pm2$ of the
supporting interval of the semicircle law, we may get improved bounds
in comparison with~(\ref{for8.1}). Namely, by Proposition \ref{pro6.1} (and repeating the
argument from the proof of the inequality of Corollary \ref{cor6.2}), for all $t>0$,
\begin{eqnarray*}
\mathbf{E}\mathrm{e}^{t |F_n(x) - G(x)|} & \leq&
\mathbf{E}\mathrm{e}^{t (F_n(x) - G(x))} + \mathbf{E}\mathrm{e}^{-t (F_n(x) - G(x))} \\
& \leq&
\mathrm{e}^{t (F(x+h) - G(x))} + \mathrm{e}^{t (G(x) - F(x-h))} \\
& \leq&
2 \mathrm{e}^{t (G(x+h) - G(x-h)) + t \|F-G\|},
\end{eqnarray*}
where $h = \sqrt{\frac{2\sigma_n^2 t}{n}} = \frac{2\sigma\sqrt{t}}{n}$.
Taking the logarithm and applying Jensen's inequality, we arrive at
%
\begin{equation}\label{for8.2}
\mathbf{E}|F_n(x) - G(x)| \leq
\|F-G\| + \bigl(G(x+h) - G(x-h)\bigr) + \frac{\log2}{t}.
\end{equation}
Using the Lipschitz property of $G$ only (that is,
$G(x+h) - G(x-h) \leq\frac{h}{\curpi}$) would yield the previous bound,
such as the one in the estimate (\ref{for6.7}) of Proposition \ref{pro6.4},
%
\begin{equation}\label{for8.3}
\mathbf{E}|F_n(x) - G(x)| \leq\|F-G\| +
C \biggl(\frac{\sigma}{n}\biggr)^{ 2/3}.
\end{equation}
However, the real size of increments $G(x+h) - G(x-h)$ with respect to
the parameter $h$ essentially depends on the point $x$. To be more careful
in the analysis of the right-hand side of~(\ref{for8.2}), we may use the following
elementary calculus bound, whose proof we omit.

\begin{lemma}\label{lem8.1} $G(x+h) - G(x-h) \leq2 g(x) h + \frac
{4}{3\curpi} h^{3/2}$  for all $x \in\mathbf{R}$ and $h>0$.
\end{lemma}

Since $G$ is concentrated on the interval $[-2,2]$, for $|x| \geq2$,
we have a simple bound $G(x+h) - G(x-h) \leq\frac{4}{3\curpi} h^{3/2}$.
As a result, one may derive from (\ref{for8.2}) an improved variant of (\ref{for8.3}).
In particular, if $|x| \geq2$, then
\[
\mathbf{E}|F_n(x) - G(x)| \leq\|F-G\| + C \biggl(\frac{\sigma}{n}\biggr)^{ 6/7}.
\]
The more general statement for all $x \in\mathbf{R}$ is given by the
following result.

\renewcommand{\thetheorem}{\arabic{section}.\arabic{theorem}}
\setcounter{theorem}{1}
\begin{theorem}\label{teo8.2}  Let
$\xi_{jk}$ $(1 \leq j \leq k \leq n)$
be independent and satisfy a logarithmic Sobolev inequality with
constant $\sigma^2$, with $\mathbf{E}\xi_{jk} = 0$ and
$\operatorname{Var}(\xi_{jk}) = 1$.
For all $x \in\mathbf{R}$,
%
\begin{equation}\label{for8.4}
\mathbf{E}|F_n(x) - G(x)| \leq\|F-G\| + C\biggl[
\biggl(\frac{\sigma}{n}\biggr)^{ 6/7} +
g(x)^{2/3}\biggl(\frac{\sigma}{n}\biggr)^{ 2/3}\biggr],
\end{equation}
where $C$ is an absolute constant.
\end{theorem}

A similar uniform bound may also be shown to hold for
$\mathbf{E}\sup_{y \leq x} |F_n(y) - G(y)|$ $(x \leq0$) and
$\mathbf{E}\sup_{y \geq x} |F_n(y) - G(y)|$ $(x \geq0$).
Note that in comparison with (\ref{for8.3}), there is an improvement for the
points $x$ at distance not more than
$(\frac{\sigma}{n})^{4/7}$ from $\pm2$.

\begin{pf*}{Proof of Theorem \ref{teo8.2}} According to the bound (\ref{for8.2}) and Lemma \ref{lem8.1}, for any $h>0$,
we may write $\mathbf{E}|F_n(x) - G(x)| \leq\|F-G\| + 3 \varphi(h)$,
where
$\varphi(h) = g(x) h + h^{3/2} + \frac{\varepsilon }{h^2}$,
$\varepsilon = (\frac{\sigma}{n})^{2}$.

We shall now estimate the minimum of this function. Write
$h = (\frac{\varepsilon }{1 + \alpha})^{2/7}$ with parameter $\alpha>0$
to be specified later on. If $g(x) \leq\alpha\sqrt{h}$, then
%
\begin{equation}\label{for8.5}
\varphi(h) \leq(1 + \alpha) h^{3/2} + \frac{\varepsilon }{h^2} =
2 (1 + \alpha)^{4/7} \varepsilon ^{3/7}.
\end{equation}
Note that the requirement on $g(x)$ is equivalent to
$\frac{g(x)^7}{\varepsilon } \leq\frac{\alpha^7}{1 + \alpha}$.
Thus, we set
$A = \frac{g(x)^7}{\varepsilon }$ and take $\alpha= 1 + 2A^{1/6}$. Since
$\alpha\geq1$, we get
$\frac{\alpha^7}{1 + \alpha} \geq\frac{\alpha^6}{2} \geq A$.
Hence, we may apply (\ref{for8.5}). Using $(1 + \alpha)^{4/7} \leq(2 \alpha)^{4/7}$
and $\alpha^{4/7} \leq1 + ((2A)^{1/6})^{4/7} = 1 + (2A)^{2/21}$,
we finally get that
\[
\varphi(h) \leq2 \cdot2^{4/7} \bigl(1 + (2A)^{2/21}\bigr) \varepsilon ^{3/7}
\leq
4 (\varepsilon ^{3/7} + A^{2/21} \varepsilon ^{3/7}).
\]
This is the desired expression in square brackets in (\ref{for8.4}) and
Theorem \ref{teo8.2} follows.
\end{pf*}

Finally, let us comment on the meaning of the general Corollary \ref{cor6.5}
in the matrix model above. To every interval $I$ on the real line,
we associate the number
\[
N_I = \operatorname{card}\{i \leq n\dvtx X_i \in I\}
\]
of eigenvalues $X_i$ inside it. Again, Corollary \ref{cor6.5} may be applied
to the standard semicircle law~$G$ with density $g(x)$, in which case
$\beta= C (\frac{\sigma}{n})^{2/3}$. This gives that, under
$\operatorname{LSI}(\sigma^2)$, imposed on the entries $\xi_{jk}$,
for any $\delta> 0$ and any interval $I$
of length $|I| \geq4 \|F-G\|/\delta$, we have that
%
\begin{equation}\label{for8.6}
\biggl|\frac{N_I}{n} - \int_I g(x) \,\mathrm{d}x \biggr| \leq\delta|I|
\end{equation}
with probability at least
$1 - 4 \exp\{ -\frac{c n^2 \delta^3}{\sigma^2} |I|^{3}\}$.
As we have already mentioned, under $\operatorname{PI}(\sigma^2)$, one can show
that $\|F-G\| \leq Cn^{-2/3}$ (\cite{B-G-T}, Theorem~\ref{teo1.1}). Therefore, (\ref{for8.6})
holds true with high probability, provided that
%
\begin{equation}\label{for8.7}
|I| \geq Cn^{-2/3}/\delta
\end{equation}
with large $C$ (of order, say, $\log^\varepsilon n$).

Such properties have been intensively studied in recent years
in connection with the universality problem. In particular, it is shown
in \cite{E-S-Y} and \cite{T-V} that the restriction (\ref{for8.7}) may be weakened to
$|I| \geq C_\varepsilon (\log^\alpha n)/n$ under the assumption that the
intervals $I$ are contained in $[-2-\varepsilon, 2+\varepsilon ]$,
$\varepsilon >0$, that is,
`in the bulk'.
\begin{appendix}
\section*{Appendix}\label{app}
\setcounter{equation}{0}

Here we recall some facts about Poincar\'e-type and log-Sobolev
inequalities. While Lemmas~\ref{lem2.2} and \ref{lem5.3} list some of their consequences,
one might wonder which measures actually satisfy these analytic
inequalities. Many
interesting examples can be constructed with the help of the following
elementary proposition.

\renewcommand{\theproposition}{A\arabic{proposition}}
\setcounter{proposition}{0}
\begin{proposition}\label{pra1} Let $\mu_1,\ldots,\mu_N$ be
probability measures
on $\mathbf{R}$ satisfying $\operatorname{PI}(\sigma^2)$ (resp., $\operatorname{LSI}(\sigma^2)$).
The image $\mu$ of the product measure $\mu_1 \otimes\cdots\otimes
\mu_N$
under any map $T\dvtx\mathbf{R}^N \rightarrow\mathbf{R}^n$ with finite
Lipschitz seminorm
satisfies $\operatorname{PI}(\sigma^2\|T\|_{\mathrm{Lip}}^2)$ (resp.,
$\operatorname{LSI}(\sigma^2\|T\|_{\mathrm{Lip}}^2)$).
\end{proposition}

On the real line, disregarding the problem of optimal constants,
Poincar\'e-type inequalities may be reduced to Hardy-type inequalities
with weights. Necessary and sufficient conditions for a measure on the
positive half-axis to satisfy a Hardy-type inequality with general weights
were found in the late 1950s in the work of Kac and Krein \cite{K-K}.
We refer the interested reader to \cite{Mu} and \cite{Maz} for a full characterization
and an account of the history; here, we just recall the principal result
(see also \cite{B-G2}).

Let $\mu$ be a probability measure on the line with median $m$, that is,
$\mu(-\infty,m) \leq\frac{1}{2}$ and $\mu(m,+\infty) \leq\frac{1}{2}$.
Define the quantities
\begin{eqnarray*}
A_0(\mu) &=& \sup_{x<m} \biggl[
\mu(-\infty,x) \int_{-\infty}^x \frac{\mathrm{d}t}{p_\mu(t)} \biggr],
\\
A_1(\mu) &=& \sup_{x>m} \biggl[
\mu(x,+\infty) \int_x^{+\infty} \frac{\mathrm{d}t}{p_\mu(t)} \biggr],
\end{eqnarray*}
where $p_\mu$ denotes the density of the absolutely continuous
component of $\mu$ (with respect to Lebesgue measure) and where
we set $A_0 = 0$ (resp., $A_1 = 0$) if
$\mu(-\infty,m) = 0$ (resp., $\mu(m,+\infty) = 0$).
We then have the following proposition.

\begin{proposition}\label{proA2} The measure $\mu$ on $\mathbf{R}$ satisfies
$\operatorname{PI(\sigma^2)}$ with some finite constant if and only
if both $A_0(\mu)$ and $A_1(\mu)$ are finite. Moreover,
the optimal value of $\sigma^2$ satisfies
\[
c_0\bigl(A_0(\mu) + A_1(\mu)\bigr) \leq\sigma^2 \leq c_1\bigl(A_0(\mu) + A_1(\mu
)\bigr),
\]
where $c_0$ and $c_1$ are positive universal constants.
\end{proposition}

Necessarily, $\mu$ must have a non-trivial absolutely continuous
part with density which is positive almost everywhere on the
supporting interval.

For example, the two-sided exponential measure $\mu_0$, with density
$\frac{1}{2} \mathrm{e}^{-|x|}$, satisfies PI($\sigma^2$) with $\sigma^2 = 4$.
Therefore, any Lipschitz transform $\mu= \mu_0 T^{-1}$ of $\mu_0$
satisfies PI($\sigma^2$) with $\sigma^2 = 4 \|T\|_{\mathrm{Lip}}^2$.\vspace*{2pt}
The latter property may be expressed analytically in terms of
the reciprocal to the so-called isoperimetric constant,
\[
H(\mu) = \mathop{\operatorname{ess\,inf}}_{x}
\frac{p_\mu(x)}{\min\{F_\mu(x),1 - F_\mu(x)\}},
\]
where $F_\mu(x) = \mu(-\infty,x]$ denotes the distribution
function of $\mu$ and $p_\mu$ the density of its absolutely
continuous component. Namely, as a variant of the Mazya--Cheeger
theorem, we have that PI($\sigma^2)$ is valid with
$\sigma^2 = 4 /H(\mu)^2$; see \cite{B-H}, Theorem 1.3.

To roughly describe the class of measures in the case, where $\mu$ is
absolutely continuous and has a positive, continuous well-behaving
density, one may note that $H(\mu)$ and the Poincar\'e constant are
finite, provided that the measure has a finite exponential moment. In
particular, any probability measure with a logarithmically concave
density satisfies PI($\sigma^2)$ with a finite $\sigma$; see~\cite{Bo1}.

As for logarithmic Sobolev inequalities, we have a similar picture,
where the standard Gaussian measure represents a basic example and
plays a similar role as the two-sided exponential distribution for
Poincar\'e-type inequalities. A full description on the real line,
resembling Proposition \ref{proA2}, was given in \cite{B-G1}. Namely, for
one-dimensional probability measure $\mu$, with previous notation, we
define the quantities
\begin{eqnarray*}
B_0(\mu) & = & \sup_{x<m} \biggl[
\mu(-\infty,x) \log\frac{1}{\mu(-\infty,x)}
\int_{-\infty}^x \frac{\mathrm{d}t}{p_\mu(t)} \biggr],  \\
B_1(\mu) & = & \sup_{x>m} \biggl[
\mu(x,+\infty) \log\frac{1}{\mu(x,+\infty)}
\int_x^{+\infty} \frac{\mathrm{d}t}{p_\mu(t)} \biggr].
\end{eqnarray*}
We then have the following proposition.

\begin{proposition}\label{proA3} The measure $\mu$ on $\mathbf{R}$ satisfies
$\operatorname{LSI(\sigma^2)}$ with some finite constant if and only
if $B_0(\mu)$ and $B_1(\mu)$ are finite. Moreover,
the optimal value of $\sigma^2$ satisfies
\[
c_0\bigl(B_0(\mu) + B_1(\mu)\bigr) \leq\sigma^2 \leq c_1\bigl(B_0(\mu) + B_1(\mu
)\bigr),
\]
where $c_0$ and $c_1$ are positive universal constants.
\end{proposition}

In particular, if $\mu$ has a log-concave density, then
$\operatorname{LSI(\sigma^2)}$
is satisfied with some finite constant if and only if $\mu$ has
sub-Gaussian tails.
\end{appendix}
\section*{Acknowledgements}

This research was supported by NSF Grants DMS-07-06866 and SFB 701.

\printhistory


\begin{thebibliography}{10}

\bibitem{Ba}
Ball, K. (1988).
Logarithmically concave functions and sections of convex sets in
$R^n.$ \textit{Studia Math.} \textbf{88} 69--84.
\MR{0932007}

\bibitem{Bh}
Bhatria, R. (1997). \textit{Matrix Analysis. Graduate Texts in Mathematics} \textbf{169}. New York: Springer.
\MR{1477662}

\bibitem{Bo2}
Bobkov, S.G. (1999). Remarks on the Gromov--Milman inequality.
\textit{Vestn. Syktyvkar. Univ. Ser. 1 Mat. Mekh. Inform.}
\textbf{3} 15--22.
\MR{1716649}

\bibitem{Bo1}
Bobkov, S.G. (1999). Isoperimetric and analytic inequalities for log-concave probability
measures. \textit{Ann. Probab.} \textbf{27} 1903--1921.
\MR{1742893}

\bibitem{B-G-L}
Bobkov, S.G., Gentil, I. and Ledoux, M. (2001). Hypercontractivity of
Hamilton--Jacobi equations. \textit{J.~Math. Pures Appl.} \textbf{88} 669--696.
\MR{1846020}

\bibitem{B-G1}
Bobkov, S.G. and G\"otze, F. (1999). Exponential integrability and transportation
cost related to logarithmic Sobolev inequalities.
\textit{J. Funct. Anal.} \textbf{1} 1--28.
\MR{1682772}

\bibitem{B-G2}
Bobkov, S.G. and G\"otze, F. (2008). Hardy-type inequalities via Riccati and
Sturm--Liouville equations. In \textit{Sobolev Spaces in Mathematics, I}
(V. Maz'ya, ed.). \textit{Intern. Math. Series}
\textbf{8} 69--86. New York: Springer.
\MR{2508839}

\bibitem{B-G-T}
Bobkov, S.G., G\"otze, F. and Tikhomirov, A.N. (2010). On concentration of empirical
measures and convergence to the semi-circle law.
Bielefeld University. Preprint.
\textit{J. Theor. Probab.} To appear.

\bibitem{B-H}
Bobkov, S.G. and Houdr\'e, C. (1997). Isoperimetric constants for product probability
measures. \textit{Ann. Probab.} \textbf{25} 184--205.
\MR{1428505}

\bibitem{B-L1}
Bobkov, S.G. and Ledoux, M. (1997). Poincare's inequalities and Talagrand's
concentration phenomenon for the exponential distribution.
\textit{Probab. Theory Related Fields} \textbf{107} 383--400.
\MR{1440138}

\bibitem{B-L2}
Bobkov, S.G. and Ledoux, M. (2009). Weighted Poincar\'e-type inequalities for Cauchy
and other convex measures. \textit{Ann. Probab.} \textbf{37} 403--427.
\MR{2510011}

\bibitem{B-U}
Borovkov, A.A. and Utev, S.A. (1983). On an inequality and a characterization of
the normal distribution connected with it. \textit{Probab. Theory
Appl.} \textbf{28} 209--218.
\MR{0700206}

\bibitem{C-S-S}
Cambanis, S., Simons, G. and Stout, W. (1976). Inequalities for $E k(X,Y)$ when the
marginals are fixed. \textit{Z.~Wahrsch. Verw. Gebiete}
\textbf{36} 285--294.
\MR{0420778}

\bibitem{C-B}
Chatterjee, S. and Bose, A. (2004). A new method for bounding rates of convergence
of empirical spectral distributions. \textit{J. Theoret. Probab.}
\textbf{17} 1003--1019.
\MR{2105745}

\bibitem{D-S}
Davidson, K.R. and Szarek, S.J. (2001). Local operator theory, random matrices and
Banach spaces. In \textit{Handbook of the Geometry of Banach Spaces}
\textbf{I} 317--366. Amsterdam: North-Holland.
\MR{1863696}

\bibitem{D}
Dudley, R.M. (1989). \textit{Real Analysis and Probability}.
 Pacific Grove, CA: Wadsworth \& Brooks/Cole Advanced Books \&
Software.
\MR{0982264}

\bibitem{D-K-W}
Dvoretzky, A., Kiefer, J. and Wolfowitz, J. (1956). Asymptotic minimax
character of the sample distribution function and of the
classical multinomial estimator. \textit{Ann. Math. Statist.}
\textbf{27} 642--669.
\MR{0083864}

\bibitem{E-S-Y}
Erd\"os, L., Schlein, B. and Yau, H.-T. (2009). Local semicircle law and complete
delocalization for Wigner random matrices.
\textit{Comm. Math. Phys.} \textbf{287} 641--655.
\MR{2481753}

\bibitem{Ev}
Evans, L.C. (1997). \textit{Partial Differential Equations}. \textit{Graduate Studies in Math.}
\textbf{19}. Providence, RI: Amer. Math. Soc.
\MR{1625845}

\bibitem{G-T1}
G\"otze, F. and Tikhomirov, A.N. (2003). Rate of convergence to the semi-circular
law. \textit{Probab. Theory Related Fields} \textbf{127} 228--276.
\MR{2013983}

\bibitem{G-T2}
G\"otze, F. and Tikhomirov, A.N. (2005). The rate of convergence for spectra of GUE
and LUE matrix ensembles. \textit{Cent. Eur. J. Math.} \textbf{3} 666--704 (electronic).
\MR{2171668}

\bibitem{G-T-T}
G\"otze, F., Tikhomirov, A.N. and Timushev, D.A. (2007). Rate of convergence to the
semi-circle law for the deformed Gaussian unitary ensemble.
\textit{Cent. Eur. J. Math.} \textbf{5} 305--334 (electronic).
\MR{2300275}

\bibitem{G-M}
Gromov, M. and Milman, V.D. (1983). A topological application of the isoperimetric
inequality. \textit{Amer. J. Math.} \textbf{105} 843--854.
\MR{0708367}

\bibitem{G-Z}
Guionnet, A. and Zeitouni, O. (2000). Concentration of the spectral measure for large
matrices. \textit{Electron. Comm. Probab.} \textbf{5} 119--136.
\MR{1781846}

\bibitem{Gu}
Gustavsson, J. (2005). Gaussian fluctuations of eigenvalues in the GUE.
\textit{Ann. Inst. H. Poincar\'e Probab. Statist.}
\textbf{41} 151--178.
\MR{2124079}

\bibitem{H}
Hensley, D. (1980). Slicing convex bodies -- bounds for slice area in terms
of the body's covariance. \textit{Proc. Amer. Math. Soc.}
\textbf{79} 619--625.
\MR{0572315}

\bibitem{K-K}
Kac, I.S. and Krein, M.G. (1958). Criteria for the discreteness of the spectrum
of a singular string. \textit{Izv. Vys\v s. U\v cebn.
Zaved. Matematika} \textbf{2} 136--153.
\MR{0139804}

\bibitem{K}
Kim, T.Y. (1999). On tail probabilities of Kolmogorov--Smirnov statistics
based on uniform mixing processes. \textit{Statist. Probab. Lett.}
\textbf{43} 217--223.
\MR{1708089}

\bibitem{L1}
Ledoux, M. (1999). Concentration of measure and logarithmic Sobolev inequalities.
In \textit{S\'eminaire de Probabilit\'es XXXIII}. \textit{Lecture Notes in Math.}
\textbf{1709} 120--216. Berlin: Springer.
\MR{1767995}

\bibitem{L2}
Ledoux, M. (2001). \textit{The Concentration of Measure Phenomenon}.
\textit{Math. Surveys and Monographs} \textbf{89}. Providence, RI: Amer. Math. Soc.
\MR{1849347}


\bibitem{L3}
Ledoux, M. (2007). Deviation inequalities on largest eigenvalues.
In \textit{Geom. Aspects of Funct. Anal., Israel Seminar 2004--2005}.
\textit{Lecture Notes in Math.} \textbf{1910} 167--219. Berlin: Springer.
\MR{2349607}

\bibitem{Mas}
Massart, P. (1990). The tight constant in the Dvoretzky--Kiefer--Wolfowitz inequality.
\textit{Ann. Probab.} \textbf{18} 1269--1283.
\MR{1062069}

\bibitem{Mau}
Maurey, B. (1991). Some deviation inequalities. \textit{Geom. Funct. Anal.}
\textbf{1} 188--197.
\MR{1097258}

\bibitem{Maz}
Maz'ya, V.G. (1985). \textit{Sobolev Spaces}. Berlin: Springer.
\MR{0817985}

\bibitem{Mu}
Muckenhoupt, B. (1972). Hardy's inequality with weights. \textit{Studia Math.}
\textbf{XLIV} 31--38.
\MR{0311856}

\bibitem{P}
Pastur, L.A. (1973). Spectra of random selfadjoint operators.
\textit{Uspehi Mat. Nauk} \textbf{28} 3--64.
\MR{0406251}

\bibitem{R}
Ruschendorf, L. (1985). The Wasserstein distance and approximation theorems.
\textit{Z. Wahrsch. Verw. Gebiete} \textbf{70} 117--129.
\MR{0795791}

\bibitem{Se}
Sen, P.K. (1974). Weak convergence of multidimensional empirical processes
for stationary $\varphi $-mixing processes. \textit{Ann. Probab.}
\textbf{2} 147--154.
\MR{0402845}

\bibitem{So}
Soshnikov, A. (1999). Universality at the edge of the spectrum in Wigner random
matrices.
\textit{Comm. Math. Phys.} \textbf{207} 697--733.
\MR{1727234}

\bibitem{Ta}
Talagrand, M. (1996). Majorizing measures: The generic chaining.
\textit{Ann. Probab.} \textbf{24} 1049--1103.
\MR{1411488}

\bibitem{T-V}
Tao, T. and Vu, V. (2009). Random matrices: Universality of local eigenvalue
statistics. Preprint.

\bibitem{Ti}
Timushev, D.A. (2006). On the rate of convergence in probability of the spectral
distribution function of a~random matrix.
\textit{Teor. Veroyatn. Primen.} \textbf{51} 618--622.
\MR{2325551}

\bibitem{V}
Vallander, S.S. (1973). Calculations of the Vasserstein distance between probability
distributions on the line.
\textit{Teor. Verojatn. Primen.} \textbf{18} 824--827.
\MR{0328982}

\bibitem{Y}
Yoshihara, K. (1975/76). Weak convergence of multidimensional empirical processes for
strong mixing sequences of stochastic vectors.
\textit{Z. Wahrsch. Verw. Gebiete}
\textbf{33} 133--137.
\MR{0385962}

\end{thebibliography}
\end{document}